\documentclass{mrlart2e}
\def\volno{7}
\def\yrno{2000}
\def\lpageno{26}
\usepackage{amssymb}
\newtheorem{thm}{Theorem}[section]
\newtheorem{prop}{Proposition}[section]
\newtheorem{corr}{Corollary}[section]
\newtheorem{lemma}{Lemma}[section]
\newtheorem{conj}{Conjecture}[section]

\theoremstyle{definition}
\newtheorem{defe}{Definition}[section]

\theoremstyle{remark}
\newtheorem{exam}{\bf Example}[section]

\newtheorem{note}{\bf Note}[section]
\newtheorem{ntn}{\bf Notation}[section]
\newtheorem{ov}{\bf Overview}[section]
\newtheorem{rem}{\bf Remark}[section]

\numberwithin{equation}{section}

\def\peq{\underline{\prec}}
\def\tpq{\underline{\widetilde \prec}}

\def\dnd{\text{\ does not divide\ }}
\def\dd{\text{\ divides\ }}

\def\eab{e_{\beta} \otimes e_{-\alpha}}

\def\qs{$\qed$}
\def\O{\text{Ord}}

\def\Z{\mathbb Z}
\def\C{\mathbb C}
\def\Res{\text{Res}}
\def\o{\otimes}
\def\e{\varepsilon}

\def\h{{\mathfrak h}}
\def\g{{\mathfrak g}}

\begin{document}

\title{Proof of the GGS Conjecture}
\author{Travis Schedler}
\address{ 059
Pforzheimer House Mail Center, Cambridge, MA 02138.}
\email{schedler@fas.harvard.edu}
\thanks{Revision received September 19, 2000}
\maketitle

\begin{abstract}
We prove the GGS conjecture \cite{GGS} (1993), which gives a
particularly simple explicit quantization of classical $r$-matrices
for Lie algebras $\mathfrak{gl}(n)$, in terms of a matrix $R \in
Mat_n(\C) \o Mat_n(\C)$ which satisfies the quantum Yang-Baxter equation
(QYBE) and the Hecke condition, whose quasiclassical limit is $r$.
The $r$-matrices were classified by Belavin and Drinfeld in the 1980's in
terms of combinatorial objects known as Belavin-Drinfeld triples.  We
prove this conjecture by showing that the GGS matrix coincides with
another quantization from \cite{ESS}, which is a more general
construction.  We do this by explicitly expanding the product from
\cite{ESS} using detailed combinatorial analysis in terms of
Belavin-Drinfeld triples.
\end{abstract}

\section{Introduction} \label{intro}
In the 1980's, Belavin and Drinfeld classified solutions $r$ of the
classical Yang-Baxter equation (CYBE) for simple Lie algebras
$\mathfrak g$ satisfying $0 \neq r + r^{21} \in (S^2
\mathfrak{g})^{\mathfrak{g}}$ \cite{BD}.  They proved that all such
solutions fall into finitely many continuous families and introduced
combinatorial objects to label these families, Belavin-Drinfeld
triples (see Section \ref{bd}).
In 1993, Gerstenhaber, Giaquinto, and Schack attempted to
quantize such solutions for Lie algebras $\mathfrak{sl}(n).$ As a
result, they formulated a conjecture stating that certain explicitly
given elements $R_{\text{GGS}} \in Mat_n(\mathbb C) \otimes Mat_n(\mathbb C)$
satisfy the quantum Yang-Baxter equation (QYBE) and the Hecke relation
\cite{GGS}.  Specifically, the conjecture assigns a family of such
elements to any Belavin-Drinfeld triple of type $A_{n-1}$.  This
conjecture is stated in Section \ref{ggsss}.

Recently, Etingof, Schiffmann, and the author found an explicit
quantization of all $r$-matrices from the Belavin-Drinfeld list.  They
did so by twisting the coproduct in the Drinfeld-Jimbo quantum group
$U_q(\g)$.  For $\g=\mathfrak{gl}(n)$, one can evaluate the universal
$R$-matrix of the twisted $U_q(\g)$ in the vector representation of
$U_q(\g)$.  This gives an element $R_J$ of $Mat_n (\C) \o Mat_n(\C)$
which satisfies the QYBE and the Hecke relation.  This element is presented
in Section \ref{ptp}.

In this paper I show that the elements $R_J$ and
$R_{\text{GGS}}$ from \cite{ESS} and \cite{GGS}
coincide.  This proves the GGS conjecture. This is done
by explicitly expanding the formula for $R_J$ using
combinatorial techniques involving Belavin-Drinfeld
triples.  The proof occupies all of Section 2.

\begin{rem} Note that the GGS conjecture was proved in some special cases
(the Cremmer-Gervais and orthogonal disjoint cases) by Hodges in
\cite{H2} and \cite{H}.  The GGS conjecture was proved in some
additional cases (orthogonal generalized disjoint) by the author in
\cite{S2}.  Also, the disjoint case was completed in \cite{S2} by
Pavel Etingof and the author.
\end{rem}

\begin{rem}
The author actually found the matrix $R_J$ with the help of a computer
before the general twist given in \cite{ESS} was found.  The matrix
was constructed to coincide with the GGS matrix in many cases, and
motivated in part the general construction given in \cite{ESS}.  See
\cite{S} for details.  Also, many steps of this proof were motivated by
and checked with computer programs.
\end{rem}

\subsection{Belavin-Drinfeld triples} \label{bd}

Let $(e_i), 1 \leq i \leq n,$ be the standard basis for $\mathbb C^n$.
Let $\Gamma = \{e_i - e_{i+1}: 1 \leq i \leq n-1\}$ be the set of
simple roots of $\mathfrak{sl}(n)$.  We will use the
notation $\alpha_i \equiv e_i - e_{i+1}$.  Let $( , )$ denote the
inner product on $\mathbb C^n$ having $(e_i)$ as an orthonormal basis.

\begin{defe} \cite{BD}
A {\sl Belavin-Drinfeld triple of type $A_{n-1}$} is a
triple\linebreak
$(T, \Gamma_1, \Gamma_2)$ where $\Gamma_1, \Gamma_2 \subset \Gamma$
and $T: \Gamma_1 \rightarrow \Gamma_2$ is a bijection, satisfying
two relations:

(a) $T$ preserves the inner product: $\forall \alpha, \beta \in
\Gamma_1$, $(T \alpha,T \beta) = (\alpha, \beta)$.

(b) $T$ is nilpotent: $\forall \alpha \in \Gamma_1, \exists k
\in \mathbb N$ such that $T^k \alpha \notin \Gamma_1$.
\end{defe}

Let $\mathfrak g = \mathfrak{gl}(n)$ be the Lie algebra of complex $n
\times n$ matrices. Let $\mathfrak h \subset \mathfrak g$ be the
subspace of diagonal matrices.  Elements of $\mathbb C^n$ define
linear functions on $\mathfrak h$ by $\bigl ( \sum_i \lambda_i e_i
\bigr) \bigl( \sum_i a_i\: e_{ii} \bigr)= \sum_i \lambda_i a_i$.  Let
$P = \sum_{1 \leq i,j \leq n} e_{ij} \otimes e_{ji}$ be the Casimir
element inverse to the standard form on $\mathfrak g$.  It is easy to see
that $P (w \o v) = v \o w$, for any $v,w \in \C^n$. Let
$P^0=\sum_i e_{ii} \o e_{ii}$ be the projection of $P$ to $\mathfrak h
\otimes \mathfrak h$.

For any Belavin-Drinfeld triple, consider the following equations for
$s \in \mathfrak h \wedge \mathfrak h$:
\begin{gather}
\label{tr02} \forall \alpha \in \Gamma_1,
\bigl[(\alpha - T \alpha) \otimes 1 \bigr] s = \frac{1}{2}
\bigl[(\alpha + T \alpha) \otimes 1\bigr] P^0.
\end{gather}
Belavin and Drinfeld showed that solutions $r \in \mathfrak{g} \o
\mathfrak g$ of the CYBE satisfying $r + r^{21} = P$, up to
isomorphism, are given by a discrete datum (the Belavin-Drinfeld
triple) and a continuous datum (a solution $s \in \h \wedge \h$ of
\eqref{tr02}).  We now describe this classification.

For $\alpha = e_i - e_j$, set $e_\alpha \equiv e_{ij}$.  Define
$|\alpha| = |j - i|$.  For any $Y \subset \Gamma$, set $\tilde Y =
\{\alpha \in \text{Span}(Y) \mid \alpha = e_i - e_j, i < j\}$ (the set
of positive roots of the subalgebra of $\mathfrak{sl}(n)$ having $Y$
as the set of simple roots). In particular we will often use the
notation $\tilde \Gamma, \tilde \Gamma_1, \tilde \Gamma_2$.  We extend
$T$ additively to a map $\tilde \Gamma_1 \rightarrow \tilde \Gamma_2$,
i.e.  $T(\alpha+\beta)=T \alpha +T \beta$. Whenever $T^k \alpha =
\beta$ for $k \geq 1$, we say $\alpha \prec \beta$.  Clearly $\prec$
is a partial ordering on $\tilde \Gamma$.  We will also use $\alpha
\peq \beta$ to denote $\alpha \prec \beta$ or $\alpha = \beta$.
Suppose $T^k \alpha = \beta$ for $\alpha = e_i - e_j$ and $\beta = e_l
- e_m$.  Then there are two possibilities on how $T^k$ sends $\alpha$
to $\beta$, since $T$ is an automorphism of the Dynkin
diagram. Namely, either $T^k(\alpha_i) = \alpha_l$ and
$T^k(\alpha_{j-1})=\alpha_{m-1}$, or $T^k(\alpha_i) = \alpha_{m-1}$ and
$T^k(\alpha_{j-1}) = \alpha_l$.  In the former case, call $T^k$ {\sl
orientation-preserving on $\alpha$}, and in the latter, {\sl
orientation-reversing on $\alpha$}.  Let
\begin{equation}
C_{\alpha,\beta} =
\begin{cases} 1, & \text{if $T^k$ reverses orientation on $\alpha$,}
\\ 0, & \text{if $T^k$ preserves orientation on $\alpha$.}
\end{cases}
\end{equation}
Now we define
\begin{gather}
a = \sum_{\alpha \prec \beta} (-1)^{C_{\alpha,\beta} (|\alpha|-1)}
(e_{-\alpha} \o e_\beta - e_{\beta} \o e_{-\alpha}),
\\ r_{st} = \frac{1}{2} \sum_i e_{ii} \o e_{ii} +
\sum_{\alpha \in \tilde \Gamma} e_{-\alpha} \o e_{\alpha}, \quad r = s + a +
r_{st}, \label{r}
\end{gather}
($r_{st} \in \mathfrak{g} \o \mathfrak{g}$ is the standard solution
of the CYBE satisfying $r_{st} + r_{st}^{21} = P$.)  The element $r$ is the
solution of the CYBE corresponding to the data
$((\Gamma_1,\Gamma_2,T), s)$.  It follows from \cite{BD} that any
solution $\tilde r \in \mathfrak{g}, \tilde r+\tilde r^{21} = P$ is
equivalent to such a solution $r$ under an automorphism of
$\mathfrak{g}$.

\subsection{The GGS conjecture} \label{ggsss}
The GGS conjecture suggests a quantization of the matrix
$r$ given in \eqref{r}, given by a matrix $R \in Mat_n(\C) \o
Mat_n(\C)$ conjectured to satisfy the quantum Yang-Baxter equation
(QYBE), $R^{12} R^{13} R^{23} = R^{23} R^{13} R^{12}$, and the Hecke
relation, $(PR - q)(PR+q^{-1}) = 0$.  This may be formulated and justified
as follows (which is more or less the original motivation).

If we write $R \equiv 1 + 2 \hbar r + 4 \hbar^2 t \pmod{\hbar^3}$,
where $q \equiv e^\hbar$, then we can consider the constraints imposed
by the QYBE and the Hecke relation modulo $\hbar^3$.  One may easily
check that the QYBE becomes the CYBE for $r$, while the
Hecke relation becomes the condition $t + t^{21} = r^2$.
Since it is not difficult to see that $r^2$ is symmetric,
the unique symmetric choice for $t$ is
$t = \frac{1}{2} r^2 = \frac{1}{2} (s^2 + (a + r_{st}) s + s (a +
r_{st}) + \varepsilon)$ where
\begin{equation} \label{eps}
\varepsilon = ar_{st} + r_{st} a + a^2.
\end{equation}

\begin{prop} \label{ggs1}
There exist unique polynomials $P_{i,j,k,l}$ of the form \\
$x q^y (q-q^{-1})^z, x,y \in \C, z \in \{0,1\}$ such that
$\sum_{i,j,k,l} P_{i,j,k,l} e_{ij} \o e_{kl} \equiv 1 + 2\hbar r
+2 \hbar^2 r^2 \pmod{\hbar^3}$.
\end{prop}
\begin{proof}  The proof is easy.  \end{proof}

\begin{defe} Define $R_{\text{GGS}} = \sum_{i,j,k,l} P_{i,j,k,l} e_{ij} \o
e_{kl}$, with the $P_{i,j,k,l}$ uniquely determined by Proposition \ref{ggs1}.
The matrix $R_{\text{GGS}}$ is called the GGS $R$-matrix.
\end{defe}

We will use the notation $x = \sum_{i,j,k,l} x_{ik}^{jl} e_{ij} \o
e_{kl}$ for elements $x \in Mat_n(\C) \o Mat_n(\C)$.  Define the 
following matrices:

\begin{equation} \label{brg}
\tilde a = \sum_{i,j,k,l} a_{ik}^{jl} q^{a_{ik}^{jl} \varepsilon_{ik}^{jl}}
e_{ij} \o e_{kl},
\quad \bar R_{\text{GGS}} = R_{st} + (q-q^{-1}) \tilde a,
\end{equation}
where
$R_{st} = q \sum_{i} e_{ii} \otimes e_{ii} +
\sum_{i \neq j} e_{ii} \otimes e_{jj} + (q - q^{-1}) \sum_{i>j} e_{ij}
\otimes e_{ji}$ is the standard Drinfeld-Jimbo solution to the QYBE,
which is a quantization of $r_{st}$.

\begin{prop}  The matrix $R_{\text{GGS}}$ equals $q^{s} \bar R_{\text{GGS}}
\, q^{s}$.
\end{prop}

\begin{proof} This is a straightforward computation. \end{proof}

\begin{rem} We see that $R_{\text{GGS}} \equiv q^{2r} \pmod{\hbar^3}$,
although $R_{\text{GGS}} \neq q^{2r}$ in general.
\end{rem}

\begin{conj}{\bf ``the GGS conjecture'' \cite{GGS}} \label{ggs}
The matrix $R_{\text{GGS}}$ satisfies the QYBE and the Hecke relation.
\end{conj}

\begin{rem} \label{or0} It is sufficient to check the QYBE for one value
of $s$ since the space of solutions to the homogeneous equation
corresponding to \eqref{tr02} is exactly the space $\Lambda^2
\mathfrak{l}$ where $\mathfrak{l} \subset \mathfrak{h}$ is the space
of all $x$ such that $(x,\alpha) = (x, T \alpha)$ for any $\alpha \in
\Gamma_1$.  Indeed, it is easy to see that $x \in \mathfrak{l}$
implies $[1 \o x + x \o 1, R_{\text{GGS}}] = 0$, and it follows that
$q^y R_{\text{GGS}} q^y$
satisfies the QYBE iff $R_{\text{GGS}}$ does, for any $y \in 
\Lambda^2 \mathfrak{l}$.
\end{rem}

\begin{rem}  Our formulation is from \cite{GH}, correcting
misprints.  The original formulation in \cite{GGS} is somewhat
different. We will write $x_{q^{-1}}$ to denote the matrix $x$ with
$q^{-1}$ substituted for $q$. Define $(x \o y)^T = x^T \o y^T$ where
$x^T$ is the transpose of $x$, for $x, y \in Mat_n(\C)$. Then, the
original form of $R_{\text{GGS}}$ can be written as follows:
\begin{equation}
R = q^{-s} \bigl( R_{st} + (q^{-1} - q) \tilde a_{q^{-1}}^T \bigr)
q^{-s}.
\end{equation}
We have $R_{\text{GGS}} - R^T_{q^{-1}} = q^{s} (q - q^{-1}) P q^{s} =
(q - q^{-1}) P$.  Thus, $R_{\text{GGS}}$ satisfies the Hecke relation
iff $R$ satisfies the Hecke relation.  In this case, we have $P
R^T_{q^{-1}} = (P R_{\text{GGS}})^{-1}$, so $R^T_{q^{-1}} =
(R_{\text{GGS}}^{-1})^{21}$, and thus $R$ satisfies the QYBE iff
$R_{\text{GGS}}$ does.  Thus, the two formulations are equivalent.
\end{rem}

\subsection{Passed $T$-pairs and a combinatorial formula for $\varepsilon$}
\label{ptp}
In this section we give a combinatorial formula for $\varepsilon$.
First let us introduce some definitions, which will be used in this
formula as well as in the proof of the main theorem.

\begin{defe}  A {\sl positive $T$-pair} is a pair $(T^k
\alpha, -\alpha)$  for $k > 0$.  We define the {\sl
order} to be
$\O(T^k \alpha, -\alpha) = k$.
The set of positive $T$-pairs is denoted $TP_+$.
\end{defe}

In Section 2, we will also define negative $T$-pairs, but so far we don't need
them.

\begin{defe}
If $\alpha = e_i - e_j \in \tilde \Gamma_1$ satisfies the property that
$T^k \alpha= e_j - e_{2j-i}$ for some $k > 0$, we say its
{\sl right-passing order} is $k$,
and denote it by $PO^r(\alpha) = k$.  If there is no such $k$, say
$PO^r(\alpha) = \infty$.  Similarly define {\sl
left-passing order}
$PO^l$.
\end{defe}

\begin{defe}
A positive $T$-pair $(\beta, -\alpha)$ is {\sl
right-passed} if\linebreak  $\O(\beta, -\alpha) >
PO^r(\alpha)$ and
$C_{\alpha,\beta} = C_{\alpha,PO^r(\alpha)}$.  We say that $(\beta,
-\alpha)$ is {\sl half right-passed} if $\O(\beta, -\alpha) =
PO^r(\alpha)$.  Similarly define the left versions, and we denote the
sets of all such $T$-pairs by $PTP_+^r, HPTP_+^r, PTP_+^l$, and
$HPTP_+^l$ (right-passed, half right-passed, left-passed, half
left-passed, respectively). The $+$ subscripts indicate positive
$T$-pairs.
\end{defe}

\begin{defe}
For a positive $T$-pair $(\beta, -\alpha)$, we define
$P^r(\beta, -\alpha)$ to be $1$ if $(\beta, -\alpha) \in PTP_+^r$,
$\frac{1}{2}$ if $(\beta, -\alpha) \in HPTP_+^r$, and otherwise $0$.
Similarly define $P^l(\beta, -\alpha)$.  It will be useful to define {\sl
symmetric} and {\sl anti-symmetric} versions: $P^s = P^r + P^l$, $P^a
= P^r - P^l$.
\end{defe}

This allows us to state a simple combinatorial formula for $\varepsilon$:

\begin{prop}  \label{fep} We may rewrite $\varepsilon$ as follows:
\begin{multline} \label{fe}
\varepsilon =\\ - \sum_{\alpha \prec \beta}
(-1)^{C_{\alpha,\beta} (|\alpha|-1)} \bigl[ P^s(\beta,
-\alpha) + C_{\alpha,\beta} (|\alpha|-1)
\bigr] (e_\beta \o e_{-\alpha} + e_{-\alpha} \o e_{\beta}).
\end{multline}
\end{prop}

\begin{proof} This is proved in Section \ref{es} and also follows
from the proof of the main theorem, Theorem \ref{mt} (see
Remark \ref{efr} for details). \end{proof}

\begin{corr} \label{frg} $\bar R_{\text{GGS}}$ is given as follows:
\begin{multline} \label{frge}
\bar R_{\text{GGS}} = \sum_{\alpha \prec \beta} (q -
q^{-1}) (-1)^{C_{\alpha,\beta} (|\alpha|-1)}
  \biggl[ q^{-P^s(\beta, -\alpha) - C_{\alpha,\beta}(|\alpha|-1)} 
(e_{-\alpha} \o
e_{\beta}) \\ - q^{P^s(\beta, -\alpha) +
C_{\alpha,\beta}(|\alpha|-1)} (e_{\beta} \o e_{-\alpha}) \biggr] + R_{st}
\end{multline}
\end{corr}

\begin{proof} Clear. \end{proof}

\begin{exam} \label{gcge} For a given $n$, there
are exactly $\phi(n)$ triples ($\phi$ is the Euler $\phi$-function) in
which $|\Gamma_1| + 1 = |\Gamma|$ \cite{GG}.  These are called
{\sl generalized Cremmer-Gervais} triples.  These are indexed by $m
\in \Z^+$, where $\text{gcd}(n,m) = 1$, and given by $\Gamma_1 =
\Gamma \setminus \{\alpha_{n-m}\}$, $\Gamma_2 = \Gamma \setminus
\{\alpha_m\}$, and $T(\alpha_i) = \alpha_{\text{Res}(i+m)}$, where
$\Res$ gives the residue modulo $n$ in $\{1,\ldots,n\}$.  For these
triples, there is a unique $s$ with first component having trace 0,
which is given by $s^{ii}_{ii} = 0, \forall i$, and
$s_{ij}^{ij} = \frac{1}{2} - \frac{1}{n}\text{Res}(\frac{j-i}{m})$
for $i \neq j$ (this is easy to verify directly and is also given in \cite{GG}).
With this $s$, $R_{\text{GGS}}$ has a very nice combinatorial
formula, which was conjectured by Giaquinto and checked in some cases.
We now state and prove this formula.
\begin{prop} \label{gp} $R_{\text{GGS}}$ is given as follows for
generalized
Cremmer-Gervais triples:
\begin{multline} \label{gcgr}
R_{\text{GGS}} =\\ q^{s} R_{st} q^{s} + \sum_{\alpha \prec
\beta} (q - q^{-1}) \bigl[ q^{\frac{-2\O(\beta,
-\alpha)}{n}} (e_{-\alpha} \o e_\beta) -
q^{\frac{2\O(\beta, -\alpha)}{n}} (e_\beta \o e_{-\alpha}) \bigr].
\end{multline}
\end{prop}
\begin{proof} See Appendix A. \end{proof}
\end{exam}

\subsection{The ESS twist and the main theorem} \label{j1}
In \cite{ESS}, an explicit quantization is given for any classical
$r$-matrix described in Section \ref{bd}.  This is given by a twist of
the standard coproduct on the quantum universal enveloping algebra
$U_q(\mathfrak{g})$.  In particular, in the $n$-dimensional
representation, this gives an element $J \in Mat_n(\C) \o Mat_n(\C)$
so that $R_J = q^{s} J^{-1} R_{st} J^{21} q^{s}$ satisfies the QYBE
and the Hecke relation.  In fact, $J$ is triangular, i.e. $J = 1 + N$
and $N = \sum_{\alpha, \beta \in \tilde \Gamma} N_{\alpha,\beta}
\eab$.

Suppose we are given $T$ so that $\text{max}_{x \in TP_+} \O(x) =
m$. Define the following matrices:

\begin{gather}
J_k = 1 + \sum_{\beta = T^k \alpha} (-q)^{-C_{\alpha,\beta}(|\alpha|-1)}
q^{P^a(\beta, -\alpha)} (q-q^{-1}) \eab, \\ J = \prod_{i=1}^m J_i, \quad
\bar R_J = J^{-1} R_{st} J^{21}, \quad R_J = q^s \bar R_J q^s.
\end{gather}

\begin{thm} \label{esst}
  \cite{ESS} The element $R_J$ satisfies the QYBE and the Hecke
relation.
\end{thm}

Now, we state the main theorem of this paper:

\begin{thm} \label{mt} For any given $T$ and $s$, $R_J = R_{\text{GGS}}$.
\end{thm}

This theorem clearly implies the GGS conjecture.

\section{Proof of the main theorem} \label{pmt}

\subsection{$T$-quadruples}\label{es}

\begin{ov}
In this section we will introduce combinatorial objects and arguments
which are sufficient to prove the combinatorial formula for
$\varepsilon$ \eqref{fe}.  To do this, we wish to directly expand
\eqref{eps}.  This involves expanding a quadratic expression in terms
of the form $\eab$ and $e_{-\alpha} \o e_\beta$ for $\alpha \peq
\beta$.  Most of the monomials in the expansion are zero. Thus, the
first step is to restrict our attention to those that are not.  In
this vein, we define {\sl compatible $T$-quadruples} (Definition
\ref{ctd}). To further simplify the formula for $\e$, we will have to
show that most of these monomials cancel pairwise.  This is
accomplished with bijections $\phi, \psi, \psi^l,$ and $\psi^r$
(Definitions \ref{phd}-\ref{pssd}) between the corresponding
$T$-quadruples.
\end{ov}

\begin{ntn}  For all of section 2 we will not need
to refer to the dimension of the representation (formerly $n$), so we will
reuse $n$ for other purposes.
\end{ntn}

\begin{defe} For any subset $Y \subset \Gamma$, define $\bar Y = (\tilde Y
\cup -\tilde Y)$.
\end{defe}


\begin{sloppypar}
\begin{defe} Let $TP_- = \{(-\alpha,\beta) \in - \tilde \Gamma
\times \tilde \Gamma \ \vert\  \alpha \peq \beta \}$.   These
will be called {\sl negative $T$-pairs}.
We define $\O(-\alpha, T^k \alpha) = k$, and $P^*(-\alpha, T^k \alpha)
= P^*(T^k \alpha, -\alpha)$ for $* \in \{l,r,a,s\}$.
\end{defe}
\end{sloppypar}

\begin{note} For negative $T$-pairs, we allow the order to be zero,
but not for positive ones!
\end{note}

\begin{defe} For $x \in TP$, define $S_x = 1$ if $x \in TP_+$ and $S_x = -1$
if $x \in TP_-$.
\end{defe}

\begin{defe} We define $PTP_-^*$ and $HPTP_-^*$ just as in
the positive case, but with components in $T$-pairs
permuted.  Let $PTP^* = PTP_-^*
\cup PTP_+^*$ and $HPTP^* = HPTP_-^* \cup HPTP_+^*$.
\end{defe}

\begin{note}  Note that all passed elements must have positive order,
so there is nothing new in the case of negative passed $T$-pairs.
\end{note}

\begin{defe} For $x = (\beta, -\alpha) \in TP_+$, set $x_l =
\alpha$ and $x_h = \beta$ ({\sl lower, higher,}
respectively.)  For $x = (-\alpha,
\beta) \in TP_-$, similarly set $x_l = \alpha$ and $x_h =
\beta$.
\end{defe}

\begin{lemma} An element $x \in TP$ is uniquely given by $x_l, x_h$, and
the sign of $x$.
\end{lemma}

\begin{proof} This is obvious. \end{proof}

For convenience, we will often give elements $x \in TP$ in terms of
$x_l, x_h$, and the sign of $x$.

\begin{defe} For convenience, we will say that $PO^l(x) = PO^l(x_l)$ and
$PO^r(x) = PO^r(x_l)$ for $T$-pairs $x$.
\end{defe}


\begin{defe}
For any pair $x = (\alpha,\beta) \in \bar \Gamma \times \bar \Gamma$,
define $E_x$ to be $E_x = e_\alpha \o e_\beta$.  We will, however,
only consider cases where $x = (e_i - e_j, e_k - e_l)$ and $i + k = j
+ l$.
\end{defe}

\begin{lemma} \label{ql} Suppose $x_1,\ldots,x_n \in TP$.
Then, if $E_{x_1} \cdots E_{x_n} \neq 0$, then there is a unique $z
\in \bar \Gamma \times \bar \Gamma$ such that $E_z = E_{x_1} \cdots
E_{x_n}$. Namely, this is given by componentwise addition of all $x_i$.
\end{lemma}

\begin{proof} This follows from nilpotency. It is easy to see
that we need only show that $E_{x_1} \cdots E_{x_n} \neq e_{ii} \o
e_{jj}$ for any $i$ and $j$.  Equivalently, we have to show that the
componentwise addition of all $x_i$ is not $(0,0)$.

First, we generalize $\O$ and $\peq$.  Clearly,
every $\alpha \in \tilde \Gamma$ can be written as $\alpha =
\sum_{i=1}^m \alpha_{k_i}$ for some $k_i$.  We say that $\alpha
\tpq \beta$ for $\alpha, \beta \in \tilde \Gamma$ if $\alpha =
\sum_{i=1}^m \alpha_{k_i}$ and $\beta = \sum_{i=1}^m T^{l_i}
\alpha_{k_i}$ for nonnegative integers $l_i$.  In this case, we define
$\widetilde{\O}(-\alpha, \beta) = \widetilde \O(\beta, -\alpha) =
\sum_i l_i$.  Note that, when $x \in TP$, $\widetilde \O(x) = |x|
\O(x)$.

Now, we set $\widetilde {TP}_+ = \{ (\beta, -\alpha) \in \tilde
\Gamma \times - \tilde \Gamma \mid \alpha \underline{
\widetilde \prec}{ \beta}, \alpha \neq \beta \}$ and
$\widetilde {TP}_- = \{ (-\alpha, \beta) \in - \tilde \Gamma \times
\tilde \Gamma \mid \alpha \underline{\widetilde \prec} \beta \}$.  Let
$\widetilde {TP} = \widetilde {TP}_- \cup \widetilde {TP}_+ \cup
\{(0,0)\}$ and let $\widetilde \O(0,0) = 0$.  We see that componentwise
addition of two elements of $\widetilde{TP}$, if in $\bar \Gamma
\times \bar \Gamma \cup \{(0,0)\}$, yields another element of
$\widetilde{TP}$, where the orders are summed.  In particular, this
means that, when $x+y = (0,0)$, for $x, y \in \widetilde {TP}$, then $x
= y = (0,0)$.  Since $TP \subset \widetilde {TP}$, it follows that any
componentwise sum of elements of $TP$, if in $\bar \Gamma \times \bar
\Gamma$, yields a nonzero element of $\widetilde {TP}$. This is all
we need.
\end{proof}

\begin{defe} Define $Q: TP^n \rightarrow TP \cup \{0\}$ as follows.
When $E_{x_1} \cdots E_{x_n} = E_z$, $z \in \bar \Gamma
\times \bar \Gamma$, for $x_i \in TP$,
$Q(x_1,\ldots,x_n) = z$.  Otherwise, $Q(x_1,\ldots,x_n) = 0.$
\end{defe}

We now apply these general definitions to the task of setting aside
the important subsets of $TP \times TP$.  These subsets
parameterize nonzero terms that arise upon expansion of $\e$.

\begin{defe} \label{ctd}
Define $CTQ, CTQ_s,$ and $CTQ_o$ ({\sl compatible $T$-quadruples of
same/opposite sign}) by
\begin{gather}\label{ctqd}
CTQ = \{ (x,y) \in TP \times TP \mid E_x E_y \neq 0\}, \\ CTQ_s =
CTQ \cap (TP_+ \times TP_+ \cup TP_- \times TP_-),
\\ CTQ_o = CTQ \cap (TP_+ \times TP_- \cup TP_- \times TP_+).
\end{gather}
\end{defe}

\begin{defe}
Define $CTQ_s^i, CTQ_s^c, CTQ_s^d$ ({\sl increasing, constant,} and
{\sl decreasing} in {\bf order}) by
\begin{gather}\label{ctqsd}
CTQ_s^i = \{ (x,y) \in CTQ_s \mid \O(x) < \O(y) \}, \\
CTQ_s^c = \{ (x,y) \in CTQ_s \mid \O(x) = \O(y) \}, \\
CTQ_s^d = \{ (x,y) \in CTQ_s \mid \O(x) > \O(y) \}.
\end{gather}
\end{defe}

\begin{lemma} For any $(x,y) \in CTQ_o$, $|x| \neq |y|$. \end{lemma}

\begin{proof} If $|x| = |y|$, then  it is clear from $(x,y) \in CTQ_o$
that $x = -y$, negating in each component.  But then $\O(x) = 
-\O(y)$, which is impossible.
\end{proof}

\begin{defe}
Define $CTQ_o^i$ and $CTQ_o^d$ ({\sl increasing} and
{\sl decreasing} in {\bf size}) by
\begin{gather}
CTQ_o^i = \{ (x,y) \in CTQ_o \mid |x| < |y| \}, \\
CTQ_o^d = \{ (x,y) \in CTQ_o \mid |x| > |y| \}.
\end{gather}
\end{defe}

\begin{note} $CTQ_o$ and $CTQ_s$ are partitioned into increasing and
decreasing quadruples by different criteria (size and order, respectively)!
It turns out that these criteria are the useful ones.
\end{note}

Now we define {\sl good} and {\sl bad} subsets of the sets we have defined.

\begin{defe}
\begin{gather}
GCTQ_s = \{(x,y) \in CTQ_s \mid \O(x)-\O(y) \dnd \O(x) \}, \\
GCTQ_o = \{(x,y) \in CTQ_o \mid |x| - |y| \dnd |x| \},
\end{gather}
and $BCTQ_* = CTQ_* \setminus GCTQ_*.$  Moreover, set
$GCTQ_{*_1}^{*_2} = GCTQ_{*_1} \cap CTQ_{*_1}^{*_2}$ and similarly for $BCTQ$.
\end{defe}

\begin{note} The definitions above are symmetric in the two components
of $T$-quadruples because $(a-b) \dd a$ iff $(a-b) \dd b$ iff $(b-a)
\dd a$ iff $(b-a) \dd b.$
\end{note}


\begin{note} In the following definitions, particularly in the
orientation-reversing cases, it may not be completely obvious that the
maps have the images indicated (i.e. are well-defined).  See Proposition
\ref{qc} for well-definition.
\end{note}

\begin{defe} \label{phd}
We define $\phi: CTQ_o \rightarrow CTQ_s \cup PTP \cup HPTP_-$.  We
first consider the increasing case. Take $(x,y) \in CTQ_o^i$ and
suppose $x_l = e_j - e_{j+|x|}$, $x_h = e_k - e_{k+|x|}$, $y_l =
e_{k+|x|-|y|} - e_{k+|x|}$, and $y_h = e_j - e_{j+|y|}$.
Suppose that $|x| = p (|y|-|x|) + q$ for $0 \leq q < |y| - |x|$.  If
$q \neq 0$ (i.e. $(x,y)$ is good), then set $\phi(x,y) = (u,v)$ where
$u$ and $v$ have the same sign as $y$ and are given by
$u_l = e_{k+|x|-|y|+q} - e_k$, $u_h = e_{j+|x|} - e_{j+|y|-q}$,
$v_l = e_{k+|x|-|y|} - e_{k+|x|-|y|+q}$, and $v_h = e_{j+|y|-q} - e_{j+|y|}$.
In other words, $u$ and $v$ are the unique same-sign $T$-pairs such
that $|v| = q$ and $Q(u,v) = Q(x,y)$.  In the case that $(x,y)$ is bad,
we simply set $\phi(x,y) = Q(x,y)$.

In the good decreasing case, we take $(y,x) \in GCTQ_o^d$, set $q$ and
$p$ as above, and again let $\phi(y,x) = (v,u)$ where $v,u$ are the
unique $T$-pairs with the same sign as $y$ such that $|v| = q$ and
$Q(v,u) = Q(y,x)$. In the bad decreasing case, we again set $\phi(y,x)
= Q(y,x)$.
\end{defe}

\begin{defe} We define $\psi: CTQ_s \setminus CTQ_s^c \rightarrow CTQ_o$
as follows.  Take $(u,v) \in CTQ_s^i$. If $(u,v)$ is good, i.e.
$\O(v) - \O(u)$ does not divide $\O(u)$, then write $\O(u) = p [\O(v)
- \O(u)] + q$ for $0 < q < \O(v)-\O(u)$.  In the case $(u,v)$ is bad,
i.e. $\O(v) - \O(u)$ divides $\O(u)$, we again write $\O(u) = p [\O(v)
- \O(u)] + q$, this time choosing $q = 0$ when $(u,v) \in TP_- \times
TP_-$ and $q = \O(v)-\O(u)$ when $(u,v) \in TP_+ \times TP_+$. Then,
we define $\psi(u,v) = (x,y)$ where $y$ has the same sign as $u$ and
$v$, while $x$ has the opposite sign, and $x_l = (T^{q}+T^{q +
\O(v) -
\O(u)} + \ldots+T^{q+(p-1)(\O(v)-\O(u))})(u_l+v_l)+T^{\O(u)} v_l$,
$x_h = (T^{\O(v)-\O(u)} + T^{2[\O(v) - \O(u)]} + \ldots +
T^{p(\O(v)-\O(u))})(u_l+v_l)+T^{(p+1)(\O(v)-\O(u))} v_l$, $y_l =
(1+T^{\O(v)-\O(u)}+\ldots+T^{p(\O(v)-\O(u))})(u_l +
v_l)+T^{(p+1)(\O(v)-\O(u))} v_l$, and $y_h =
(T^q+T^{q+(\O(v)-\O(u))}+\ldots+T^{\O(u)})\linebreak (u_l +
v_l)+T^{\O(v)} v_l$.

In the case of decreasing quadruples, we begin with $(v,u) \in
CTQ_s^d$ and set $\psi(v,u) = (y,x)$, with $y, x, p,$ and $q$ all defined as
above.
\end{defe}

\begin{defe} \label{pssd} We similarly define $\psi^r: PTP^r \cup 
HPTP^r_- \rightarrow
CTQ_o^i$ and $\psi^l: PTP^l \cup HPTP^l_- \rightarrow
CTQ_o^d$. Suppose $v \in PTP^r \cup HPTP^r_-$. If $v$ is good,
i.e. $PO^r(v)$ does not divide $\O(v)$, we write $\O(v) = p\, PO^r(v) +
q$ for $0 < q < PO^r(v)$.  If $v$ is bad, i.e.  $PO^r(v)$ divides
$\O(v)$, we set $q = 0$ when $v \in TP_-$ and $q = PO^r(v)$ when $v
\in TP_+$, and again write $\O(v) = p\, PO^r(v) + q$.  Now, we define
$\psi^r(v) = (x,y)$ where $y$ has the same sign as $v$, $x$ has the
opposite sign, and $x_l = (T^{q} + T^{q + PO^r(v)} + \ldots +
T^{q + (p-1) PO^r(v)}) v_l, x_h = (T^{PO^r(v)} + T^{2PO^r(v)} + \ldots +
T^{p\, PO^r(v)}) v_l, y_l = (1 + T^{PO^r(v)} + \ldots + T^{p\, PO^r(v)}) v_l,$
and $y_h = (T^{q} + T^{q+PO^r(v)} + \ldots + T^{q+p\, PO^r(v)}) v_l$.
For $v \in PTP^l \cup HPTP^l_-$, we define $\psi^l(v) = (y,x)$, with
$y, x, p,$ and $q$ all defined as above.
\end{defe}

\begin{note} We had to define two separate maps $\psi^r$ and $\psi^l$ because
sometimes $HPTP^r_- \cup PTP^r$ and $HPTP^l_- \cup PTP^l$ intersect, and
$\psi^l$ and $\psi^r$ do not agree.
\end{note}

\begin{defe}  We say that a pair $x = (\beta, -\alpha) \in TP_+$ reverses
orientation if $C_{\alpha,\beta} = 1$, and in this case, we set $C_x =
1$.  Otherwise $x$ preserves orientation and $C_x = 0$.  For negative
pairs $x = (-\alpha, \beta) \in TP_-$, we say $x$ preserves
orientation and $C_x = 1$ if $\alpha = \beta$; otherwise, we set $C_x
= C_{(\beta, -\alpha)}$ and say that $x$ preserves/reverses
orientation iff $(\beta, -\alpha)$ does.
\end{defe}

\begin{defe} Define $s_q, s: TP^n \rightarrow \Z[q,q^{-1}]$ as follows:
$s_q(x_1, \ldots, x_n) = (-q)^{-C_{x_1} (|x_1| - 1) - \ldots - 
C_{x_n} (|x_n| - 1)}$ and $s = s_1$.
\end{defe}

\begin{lemma} \label{nc} If $(u, v) \in CTQ_s^i$ then $u$ preserves 
orientation.  In the event $v$ reverses orientation, then $\O(u) < 
\frac{\O(v)}{2}$.
If $(x,y) \in CTQ_o^i$, then $y$ preserves orientation. In the event
$x$ reverses orientation, $|x| \leq \frac{|y|}{2}$.  The same results
hold considering $(v, u) \in CTQ_s^d$ and $(y, x) \in CTQ_o^d$,
respectively.
\end{lemma}

\begin{proof}  In the first case, if $u$ reversed orientation, then
$T^{\O(u)}(v_l) = v_h$, thus $\O(u) = \O(v)$ by nilpotency---this is a
contradiction.  On the other hand, when $v$ reverses orientation,
$T^{\O(v)-\O(u)}$ must reverse orientation, and thus by nilpotency
cannot be defined on all of $u_l + v_l$.  Thus $\O(v) - \O(u) <
\O(u)$, hence the desired result.  The second case follows easily from
nilpotency. \end{proof}

\begin{prop} \label{qc}
The maps $\phi \bigl|_{GCTQ_o}: GCTQ_o \rightarrow CTQ_s \setminus
CTQ_s^c$ and $\psi: CTQ_s \setminus CTQ_s^c \rightarrow GCTQ_o$ are
inverse to each other.  The maps $\phi \bigl|_{BCTQ_o^i}: BCTQ_o^i
\rightarrow PTP^r \cup HPTP^r_-$ and $\psi^r: PTP^r \cup HPTP^r_-
\rightarrow BCTQ_o^i$ are inverse to each other.  Finally, the maps
$\phi \bigl|_{BCTQ_o^d}: BCTQ_o^d \rightarrow PTP^l \cup HPTP^l_-$ and
$\psi^l: PTP^l \cup HPTP^l_- \rightarrow BCTQ_o^d$ are inverse to each
other.  All maps $\phi, \psi, \psi^r,$ and $\psi^l$ preserve $Q$ and $s_q$.
\end{prop}

\begin{proof}  The fact that the maps are well-defined,
inverse to each other, and preserve $Q$ is easy to see from
construction when orientations are preserved (it helps to draw a
picture).
Also, when orientations are preserved, $s_q$ is trivially
preserved.

So it remains to consider orientation-reversing cases.  Given $(x,y)
\in CTQ_o^i$, if orientation is reversed in $x$ or $y$, it can only be
reversed in $x$, and $|x| \leq \frac{|y|}{2}$ by Lemma \ref{nc}.  In
the case $|x| \neq \frac{|y|}{2}$, it follows that $\phi(x,y) = (u,v)
\in CTQ_s^i$, where $u$ preserves orientation, $v$ reverses
orientation, $|v| = |x|$, and $\O(y) = \O(u) < \O(v) - \O(u) = \O(x) +
\O(y)$.  So $s_q$ and $Q$ are preserved, and $\psi(u,v) = (x,y)$.  If,
instead, $|x| = \frac{|y|}{2}$, then $\phi(x,y) = v$, where $|v| =
|x|$, $v \in PTP^r \cup HPTP_-^r$, and $v$ reverses orientation.
Again, $s_q$ and $Q$ are preserved and $\psi^r(v) = (x,y)$.

Now consider an element $(u,v) \in CTQ^i_s \setminus CTQ_s^c$ in which
either $u$ or $v$ reverses orientation.  By Lemma \ref{nc}, only $v$
reverses orientation and $\O(u) < \frac{\O(v)}{2}$.  So, we get
$\psi(u,v) = (x,y) \in CTQ_o^i$ where $\O(u) = \O(y), \O(v) = 2 \O(y)
+ \O(x), |x| = |v|, |y| = 2|v| + |u|$.  In this case, $x$ reverses
orientation and $y$ does not, so $Q$ and $s_q$ are preserved, and
clearly $\phi(x,y) = (u,v)$.  Finally, suppose $v \in PTP^r \cup
HPTP^r_-$ reverses orientation.  By definition, this means that
$T^{PO^r(v)}$ reverses orientation on $v$.  By nilpotency, $\O(v) < 2
PO^r(v)$, and it follows that $\psi^r(v) = (x,y) \in CTQ_o^i$.  We
then have $|x| = |v|, \O(v) = 2\O(y) + \O(x),$ and $PO^r(v) = \O(x) +
\O(y)$.  It follows that $x$ reverses orientation, and so $s_q$ and
$Q$ are preserved, and $\phi(x,y) = v$.

The decreasing and left cases follow in exactly the same way as the
increasing and right cases.  \end{proof}

\begin{corr} \label{psq} Suppose $(x,y) \in CTQ_o^i$.  If
$\phi(x,y) = (u,v)$ and not all of $x,y,u,v$ preserve orientation,
then $x,v$ reverse orientation, $y,u$ preserve orientation, and
$|x| = |v|$. If $\phi(x,y) = v$ and not all of $v,x,y$ preserve orientation,
then $x,v$ reverse orientation, $y$ preserves orientation, and $|x| = |v|$.
The same facts hold under the assumptions $(y,x) \in CTQ_o^d$ with
$\phi(y,x) = (v,u)$ or $v$.
\end{corr}

\begin{proof} This follows directly from the argument above. \end{proof}

\begin{lemma} $\sum_{(x,y) \in CTQ_s^c} s(x,y) E_{Q(x,y)} =
\sum_{x \in TP} (1-|x|) s(x) C_x E_x.$
\end{lemma}

\begin{proof} Fix a choice of sign $\pm$ for this proof. Clearly,
whenever $E_x E_y \neq 0$, $\O(x) = \O(y)$, and $x,y \in TP_{\pm}$,
then $C_x = C_y = 1$.  In this case, $s(x) s(y) = -s(Q(x,y))$, as
$Q(x,y)$ also has reversed orientation and the same order as $x$ and
$y$.  It remains only to see that, for any $z \in TP_{\pm}$ with $C_z
= 1$, there are $|z|-1$ ways of writing $E_z = E_x E_y$ for $x,y \in
TP_{\pm}$, and they all are of this form.  The formula follows
immediately. \end{proof}

\noindent{\bf Direct proof of Proposition \ref{fep}.}  Set $b =
\sum_i e_{ii} \o e_{ii}$ and \\ $P_- = \sum_{i < j} e_{ji} \o
e_{ij}$, so that
$r_{st} = b + P_-$. Then, using $\phi$ and $\psi$,
\begin{multline}
a^2 + ar_{st} + r_{st}a = a^2 + a P_- + P_- a + \frac{1}{2} (a b + b a)
\\ = \sum_{(x,y) \in CTQ_s} s(x,y) E_{Q(x,y)} - \sum_{(x,y) \in CTQ_o} s(x,y)
E_{Q(x,y)}  - \frac{1}{2} \sum_{x \in HPTP} S_x s(x) E_x \\
= -\sum_{x \in HPTP_- \cup PTP} s(x) E_x + \sum_{(x,y) \in CTQ_s^c} 
s(x,y) E_{Q(x,y)}
- \frac{1}{2} \sum_{x \in HPTP} S_x s(x) E_x \\
= \sum_{x \in TP} (1-|x|) s(x) C_x E_x
- \sum_{x \in PTP} s(x) E_x - \frac{1}{2} \sum_{x \in HPTP} s(x) E_x.
\end{multline}

Equation \eqref{fe} follows immediately. \qs

\begin{rem} \label{efr} The combinatorial formula \ref{fe} for $\e$
also follows from the proof of Theorem \ref{mt}.  Namely, in the proof
we actually show that $R_J$ has the form of $R_{\text{GGS}}$ but we use
the combinatorial formula \eqref{fe} for $\e$ instead of the original one.
On the other hand, since the combinatorial formula \eqref{fe} is
symmetric, and because $R_J$ satisfies the QYBE and the Hecke
relation, $R_J$ must be the unique element satisfying the hypotheses
of Proposition \ref{ggs1} by the discussion in Section \ref{ggsss}.
%
The proof of Proposition \ref{fep} above is, however, given for
pedagogical reasons and because the results used will be needed later.
\end{rem}

\subsection{Passing properties of $T$-quadruples}

\begin{ov}
By {\sl passing properties} of a $T$-pair we mean information about
its left- and right-passing order.  In particular, usually we will be
concerned with whether a pair is (half) right- or (half) left-passed.

In this section, we will list all possible passing properties of
compatible increasing and decreasing quadruples, in connection with
those properties of their images under $\phi$ or $\psi$.  These results,
Lemmas \ref{pq} and \ref{ppp}, are essential in order to consider
quadratic terms which arise in the Hecke condition for
$R_{\text{GGS}}$, which are similar to those in the formula for
$\e$ but include powers of $q$ which depend on the passing
properties.  As a consequence of these results, one can prove the
Hecke condition for $R_{\text{GGS}}$ directly (see \cite{S2}).
\end{ov}

In order to prove these results, we will first need to develop some
more powerful combinatorial tools and notation regarding
Belavin-Drinfeld triples.  The combinatorics can best be pictured on
the Dynkin diagram for $\mathfrak{sl}(n)$.  We picture this diagram as
the line segment $[1,n]$ with integer vertices.  We then picture the
positive root $e_i - e_j$, for $1 \leq i < j \leq n$, as the line
segment $[i,j]$.  In this context, $\Gamma_1$ and $\Gamma_2$ can be
thought of as subsets of the graph consisting of the union of all the
length-1 segments which make them up, and maps $T$ are nilpotent
graph isomorphisms $\Gamma_1 \rightarrow \Gamma_2$.

\begin{defe} Assume $i < j$ and $k < l$.  We say $e_i - e_j < e_k - 
e_l$ if $i < k$ ($e_i - e_j$ is
to the left of $e_k - e_l$).  As subcases of this, we say that $e_i -
e_j \ll e_k - e_l$ for $j < k$, $e_i - e_j \lessdot e_k - e_l$ if $j =
k$, and $e_i - e_j \overline{<} e_k - e_l$ for $j > k$.  Similarly
define $>$ by $\alpha > \beta$ whenever $\beta < \alpha$, and the
same for $\gg, \gtrdot,$ and $\overline{>}$.
%
If $\alpha \ll \beta$ or $\beta \ll \alpha$, then we say that
$\alpha \perp \beta$
(meaning
there are orthogonal subsets $X,Y \subset \Gamma$ so that $\alpha \in
\tilde X$ and $\beta \in \tilde Y$).
\end{defe}

Here and in the sequel we will make frequent use of the following key
combinatorial lemma:

\begin{defe} Take $\alpha \in \tilde \Gamma_1$.  Let $M_\alpha$
be the smallest positive integer such that $T^{M_\alpha - 1} \alpha
\in \tilde \Gamma_2 \setminus \tilde \Gamma_1$.  Let $c_\alpha$ be the
smallest positive integer less than $M_\alpha$ such that $T^{c_\alpha}
\alpha \not \perp \alpha$, if such an integer exists. Otherwise set
$c_\alpha = \infty$.  If $c_\alpha$ is finite and there is
a positive integer $d < M_\alpha$ which is not a multiple
of $c_\alpha$ satisfying $T^d \alpha \not \perp \alpha$,
then let $d_\alpha$ be the smallest such $d$. Otherwise
set $d_\alpha = \infty$.
\end{defe}

\begin{lemma} \label{c1} Suppose $\alpha \in \tilde \Gamma_1$
and $c_\alpha, d_\alpha < \infty$.  (i) For any positive integer $d <
M_\alpha$ such that $T^d \alpha \not \perp \alpha$, either $d$ is a
multiple of $c_\alpha$, or $d = d_\alpha$. (ii) $d_\alpha + c_\alpha -
\text{gcd}(c_\alpha, d_\alpha) \geq M_\alpha$ and $T^{d_\alpha} \alpha
\perp T^{c_\alpha} \alpha.$  \end{lemma}

\begin{proof}  Take some $d$ such that $T^d \alpha \not \perp \alpha$.
We show in the following paragraph that either $c_\alpha$ divides $d$,
or $c_\alpha + d - \text{gcd} (d, c_\alpha) \geq M_\alpha$ and $T^d
\alpha \perp T^{c_\alpha} \alpha$.  This proves the lemma---all that
remains is to see that, in the latter case, $d = d_\alpha$.  If,
instead, $d \neq d_\alpha$, then applying the above result also to
$d_\alpha$, we find that both $T^d \alpha$ and $T^{d_\alpha} \alpha$
are perpendicular to $T^{c_\alpha} \alpha$.  By space concerns on the
diagram, it follows that $T^d \alpha \not \perp T^{d_\alpha} \alpha$,
but then $T^{d - d_\alpha} \alpha \not \perp \alpha$, which would show
that $c_\alpha$ divides $d - d_\alpha$, in contradiction to $d_\alpha
+ c_\alpha - \text{gcd}(c_\alpha, d_\alpha) \geq M_\alpha$.

So, take any $d$ such that $T^d \alpha \not \perp \alpha$, and assume
that $d$ is not a multiple of $c_\alpha$.  Define $f:
\{0,\ldots,{M_\alpha}-1\} \times \Gamma_1 \rightarrow \{-1,1\}$ by
$f(p, \alpha_i) = 1$ if $T^p$ preserves orientation on $\alpha_i$, and
$f(p, \alpha_i) = -1$ otherwise.  Define $g: \{0,\ldots,{M_\alpha}-1\}
\times \Gamma_1 \rightarrow \Z$ as follows. For any $1 \leq i \leq
n-1$ ($n$ is the length of the Dynkin diagram), let $q$ be given by
$T^p \alpha_i = \alpha_q$.  Then we define $g(p,\alpha_i) = q -
f(\alpha_i) i$. Define $F = f \times g$.  Clearly, $F$ is defined so
that if $Y \subset \Gamma_1$ is a connected segment of the diagram,
then $F$ is constant on $\{p\} \times Y$, for each fixed $p$, $0 \leq
p \leq M_\alpha - 1$.  Since $T^{ac_\alpha} \alpha \not \perp
T^{(a-1)c_\alpha} \alpha$ for $ac_\alpha < M_\alpha$, it follows that
$F$ is $c_\alpha$-periodic in the first component.  For the same
reason, $F$ is $d$-periodic in the first component. If
$d+c_\alpha-\text{gcd}(d, c_\alpha) < M_\alpha$, $F$ must be
$\text{gcd}(d, c_\alpha)$-periodic in the first component.  This
follows since
$F(d+l,\alpha_i)=F(l,\alpha_i)$, $0 \leq l <
c_\alpha-\text{gcd}(d,c_\alpha)$ implies $F(a,\alpha_i) =
F(b,\alpha_i)$ whenever $a \equiv b \pmod{\text{gcd}(d,c_\alpha)}$, $0
\leq a,b < {M_\alpha}$.  By minimality, $\text{gcd}(d, c_{\alpha}) =
c_\alpha$, which is impossible.  Hence, $d + c_\alpha - \text{gcd}(d,
c_\alpha) \geq M_\alpha$.  Furthermore, $T^d \perp T^{c_\alpha}$, because
otherwise $T^{d - c_\alpha} \alpha \not \perp \alpha$, which would imply
that $d - \text{gcd}(d, c_{\alpha}) \geq M_\alpha$ by the above results
applied to $d - c_\alpha$, which is clearly contradictory.
\end{proof}

\begin{corr} \label{c3}
Suppose that $T^l(\alpha_i) = \alpha_{i+r}$ for $a \leq i < b$ where $r \leq
b-a$ and $l > 0$.  Then for any $\alpha = e_{c} - e_{c+s}$, $\beta =
e_d - e_{d+s}$, where $s \geq r$ and $a \leq c < d \leq b$, then
$\beta \not \prec \alpha$, and $\alpha \prec \beta$ implies that
$\O(\beta, -\alpha) = \frac{l(d-c)}{r}$.

By reversing all directions, given $T^l(\alpha_i) = \alpha_{i-r}$, $a \leq i <
b$, with $r \leq b-a$, then for any $\alpha = e_{c-s} - e_{c}$, $\beta
= e_{d-s} - e_{d}$, where $s \geq r$ and $a \leq d < c \leq b$, then
$\beta \not \prec \alpha$, and $\alpha \prec \beta$ implies
$\O(\beta, -\alpha) = \frac{l(c-d)}{r}$.
\end{corr}

\begin{proof}  First, it is clear that $r$ divides $c-d$ iff $T^{kl} \alpha
= \beta$ for $k > 0$ an integer.  In this case, the theorem is
satisfied; so suppose not.  By applying $T^{\pm l}$ some number of
times to $\alpha$ or $\beta$, it suffices to assume $0 < d - c < r$.
Now, if $T^l$ is defined on $\alpha$, then $\alpha < \beta < T^l
\alpha$ together with the Lemma gives the desired result.  So assume
$T^l$ is not defined on $\alpha$, and hence it is not defined on
$\beta$ either.  Suppose $T^m \beta = \alpha$ for some positive
integer $m$.  Write $\beta = e_i - e_j$.  Now, by applying $T^m$ some
number of times to each $\alpha_p, i \leq p < j$, we can obtain
$\alpha_q$ for some $a \leq q \leq b$, showing that $T^l$ is defined
on $\beta$, which is a contradiction.
%
But then $T^m$ is defined
The direction-reversed case is the same. \end{proof}

\begin{defe}
For convenience, let $EPTP^* = HPTP^* \cup PTP^*$ for any or no
superscript $*$.
\end{defe}
\begin{defe}When $T^k$ acts on some segment of the Dynkin diagram
(i.e. some subset of $\Gamma$) by sending $\alpha_i$ to
$\alpha_{i+k}$, we say it acts by {\sl shifting to the right by $k$} when
$k$ is positive, and by {\sl shifting to the left by $k$} when $k$ is
negative.  In particular, on each segment, $T^k$ acts by {\sl shifting} iff
orientation is preserved.
\end{defe}

The next two lemmas summarize all of the possible passing properties
of an opposite-sign quadruple and its image under $\phi$.

\begin{lemma} \label{pq} Take $(x,y) \in GCTQ_o^i$
and
$\phi(x,y) = (u,v)$.  Then
exactly one of the following must hold:

\noindent
$($a$)$ $P^l(v) = P^r(y) + P^l(x), P^l(y) = 0$ \\
$($b$)$ $P^l(v) = P^l(y), P^l(x) = P^r(y) = 0$

Similarly, exactly one of the following must hold:

\noindent
$($c$)$ $P^r(v) = P^r(x) + P^r(u), P^l(u) = 0$ \\
$($d$)$ $P^r(v)=P^l(u), P^r(u) = P^r(x) = 0$.

These results, after interchanging superscripts of $l$ with $r$, also
  hold when one considers $(y,x) \in GCTQ_o^d$ and $\phi(y,x) = (v,u)$,
  instead of the original hypothesis.
\end{lemma}

\begin{proof}
Let $\O(v) = (k+1) \O(y) + k \, \O(x)$ and $\O(u) = k \,
\O(y) + (k-1) \O(x)$ for some $k$ (which exists by construction).

Suppose that $v \in EPTP^l$.  We will analyze all possible cases by
considering the value of $PO^l(v)$.  Write $PO^l(v) = p [\O(x) +
\O(y)] + q$ where $q < \O(x) + \O(y)$.  First, I claim that $p \geq
k-1$.  Suppose instead that $p < k-1$.  In this case, $k \geq 2$,
which immediately implies from the proof of Proposition \ref{qc} that
$T^{\O(x)+\O(y)}$ preserves orientation on $u_l + v_l$.  So $v$ preserves
orientation, and by definition $T^{PO^l(v)}$ preserves orientation on $v$.
Now, if we set $\omega = T^{-\O(y)} x_l$, we find that
$T^q \omega \overline{<} \omega$ because $p < k$. Now, $k \geq 2$
shows that $T^{\O(x)+\O(y)} \omega \overline{>} \omega$.  Finally,
$T^{\O(x)+\O(y)} \omega \overline{>} T^q \omega$, because $p < k-1$.
These facts, however, contradict Lemma \ref{c1}.

So, it must be that $p \geq k-1$.  We divide into the two cases, (1)
$p = k-1$ and (2) $p = k$.

First consider the case $k = p-1$.  Set $t = k [ \O(x) + \O(y)] -
PO^l(v)$.  Then we have three cases: (i) $t < \O(y)$, (ii) $t =
\O(y)$, and (iii) $t > \O(y)$.  First consider (i).  Now, $T^t(y_l)
\gtrdot y_l$, so $y \in PTP^r$.  Conversely, whenever $y \in PTP^r$,
clearly $v \in PTP^l$ with $PO^l(v) = k[\O(x) + \O(y)] - PO^r(y)$ (we
use that $y$ always preserves orientation).  In this case, $t < \O(y)$
and $p = k-1$, as desired.  This situation, characterized by $v
\in PTP^l, y \in PTP^r$, falls into (a) and we will call it (a1).

Next, take (ii). In this case, $PO^l(v) = k \, \O(x) + (k-1)
\O(y)$, so $x \in HPTP^l$ and $y \in HPTP^r$.
Conversely, it is clear that $y \in HPTP^r$ iff $x \in HPTP^l$ from the
construction of $\phi$, and in this case, $v \in PTP^l$ with $PO^l(v)
= k \, \O(x) + (k-1)\O(y)$.  Thus, $t = \O(y)$ and $p = k-1$, as desired.
This is a different case of (a), so let us call it (a2).

Finally, consider (iii).  In this case, $0 < \O(x) + \O(y) - t <
\O(x)$.  Set $t' = \O(x) + \O(y) - t = PO^l(v) - (k-1)[\O(x) +
\O(y)]$.  Then it follows that $T^{t'} (x_l) \lessdot x_l$, so that $x
\in PTP^l$ with $PO^l(x) = t'$.  Conversely, if $x \in PTP^l$, then $v
\in PTP^l$ with $PO^l(v) = (k-1) [\O(x) + \O(y)] + PO^l(x)$, as
desired (this can be checked separately when $x$ reverses
orientation---here $k = 1$ so there is no difficulty.)  Hence, $t >
\O(y)$ and $p = k-1$.   This is the final case of (a), so let us call it (a3).

Next, consider the case $k = p$.  Set $t = PO^l(v) - k[\O(x) +
\O(y)]$.  Because $k = p$, it follows that $0 < t \leq \O(y)$, hence 
$y \in EPTP^l$
with $PO^l(y) = t$.  Since $PO^l(y) = \O(y)$ iff $PO^l(v) = \O(v) =
(k+1)\O(y) + k \, \O(x)$, we have $P^l(v) = P^l(y)$. Conversely,
whenever $y \in EPTP^l$, then $v \in EPTP^l$ with $PO^l(v) = PO^l(y) +
k [\O(x) + \O(y)]$.  Hence, $p = k$ and $P^l(v) = P^l(y)$.  This
accounts for case (b).

We have proved the first part of the Lemma, because we have considered
all possible nonzero values of $P^l(v), P^l(x), P^l(y),$ and $P^r(y)$,
and grouped them into the cases (a1), (a2), (a3), and (b).  We have shown
that each of these is associated with different values of $PO^l(v)$,
which justifies the zero values of $P^l(v), P^l(x), P^l(y),$ and
$P^r(y)$ in each case.

Next, we apply the same analysis used in the first part to show that
exactly one of (c),(d) holds.  Let $t = | PO^r(v) - \O(x) - \O(y) |$ and
$\omega = v_l + u_l + T^{\O(x) + \O(y)} v_l$. We suppose that
$v \in EPTP^r$ and divide into the cases
$PO^r(v) > \O(x) + \O(y)$ and $PO^r(v) < \O(x) + \O(y)$.

First consider $PO^r(v) > \O(x) + \O(y)$. Then it is clear that $t
\leq \O(u)$, with equality iff $v \in HPTP^r$. Now, $T^{t}$ is defined
on $\omega$, and $T^t \omega \not \perp \omega$, so $T^t$ preserves
orientation on $\omega$.  This shows that $T^t u_l \lessdot u_l$, so
that $u \in EPTP^l$ with $PO^l(u) = t$.  That is, $P^l(u) = P^r(v)$.
Conversely, whenever $u \in EPTP^l$, we know from the fact that $u$
preserves orientation that $T^{PO^l(u)}$ preserves orientation on
$u_l$, and hence that $P^r(v) = P^l(u)$ with $PO^r(v) =
PO^l(u)+\O(x)+\O(y)$.  It is then clear that $PO^r(v) > \O(x) +
\O(y)$.  This accounts for case (d).

Now, suppose that $v \in EPTP^r$ with $PO^r(v) < \O(x) + \O(y)$.  We
divide into the cases $k \geq 2$ and $k = 1$.

First suppose $k \geq 2$. In particular, this implies that
$T^{PO^r(v)}$ is defined on $\omega$, and by nilpotency, it must
preserve orientation.  Also, $T^{\O(x)+\O(y)}$ must preserve
orientation on $\omega$.  So, we see that $T^t u_l \gtrdot u_l$, thus
$u_l \in PTP^r$ with $PO^r(u_l) = t$.  Conversely, if $k \geq 2$ and
$u \in EPTP^r$, we see from $\omega \overline{<} T^{PO^r(u)} \omega
\overline{<} T^{\O(x)+\O(y)} \omega$ and Lemma \ref{c1} that $PO^r(u)
< \O(x)+\O(y)$.  It follows that $u,v \in PTP^r$ with $PO^r(v) = t$,
and $PO^r(v) < \O(x)+\O(y)$, as desired.  This accounts for one
situation of (c); call this (c1).

Next, suppose $k = 1$ and $PO^r(v) < \O(x) + \O(y)$.  We further
divide into the three cases (i) $PO^r(v) < \O(x)$, (ii) $PO^r(v) =
\O(x)$, and (iii) $PO^r(v) > \O(x)$.

In case (i), we use that $y$ preserves orientation to see that
$T^{PO^r(v)}(x_l) \gtrdot x_l$, so that $x \in PTP^r$ with $PO^r(x) =
PO^r(v)$.  Conversely, whenever $x \in PTP^r$, it follows that $k = 1$
using Lemma \ref{c1}: otherwise we would have $T^{-\O(y)} x_l
\overline{<} x_h \overline{<} T^{PO^r(x) - \O(y)} x_l$ with $T^{-\O(y)}
x_l \lessdot T^{PO^r(x) - \O(y)} x_l$ while \\ $-\O(y) < PO^r(x) - \O(y) <
\O(x)$.  Hence, $v \in PTP^r$ with $PO^r(v) = PO^r(x) <
\O(x)$.  Call this situation (c3).


In case (ii), we have that $x, u \in HPTP^r$.  Conversely, whenever $x
\in HPTP^r$, we must have $k = 1$ by Lemma \ref{c1}, and then $u \in
HPTP^r$ and $PO^r(v) = PO^r(x)$, as desired.  If $u \in HPTP^r$, then
it also follows from Lemma \ref{c1} considering $\omega$ that $k = 1$,
and then $x \in HPTP^r$.  Call this situation (c2).

Finally, we consider case (iii).  Now, $PO^r(v) > \O(x)$ shows that
$\O(x)+\O(y)-PO^r(v) < \O(y)$.  By nilpotency,
$T^{t}$ must preserve orientation on $y_l$, and it
follows that $u \in PTP^r$ with $PO^r(u) = t$.
Conversely, if $u \in PTP^r$ with $k = 1$, then it follows that
$PO^r(u) < \O(u) = \O(y)$, so\linebreak $T^{\O(x)+\O(y)-PO^r(u)}
v_l
\gtrdot v_l$.  Hence, $v \in PTP^r$ with $PO^r(v) = \O(x)+\O(y) -
PO^r(u)$, and $\O(x) < PO^r(v) < \O(x)+\O(y)$.  This has the same
passing properties as (c1), so call this situation (c1').

We have finished the second half of the Lemma, since we have accounted
for all possible values of $P^r(u), P^l(u), P^r(v)$, and $P^r(x)$ in
cases (d), (c1), (c1'), (c2), and (c3).  Each of these are associated
with distinct values of $P^r(v)$ with rexspect to the $\O(x), \O(y),
\O(u)$, and $\O(v)$, once $k$ is fixed.

To obtain the result for decreasing quadruples, simply reverse
all directions and permute the components of all $T$-quadruples
(not pairs!) in this proof.
\end{proof}

\begin{lemma} \label{ppp}  Take $(x,y) \in BCTQ_o^i$
and $\phi(x,y) = v$.  Then $v \in EPTP^r,\linebreak P^r(v) +
P^r(x) = 1$, and exactly one of the following hold:

\noindent
$($a$)$ $P^l(v) = P^r(y) + P^l(x), P^l(y) = 0$. \\
$($b$)$ $P^l(v) = P^l(y)$, $P^l(x) = P^r(y) = 0$.

Under the hypotheses $(y,x) \in BCTQ_o^d$ and $\phi(y,x) = v$, these
results still hold upon interchanging superscripts of $l$ and $r$.
\end{lemma}

\begin{proof} Clearly $v \in HPTP^r$ iff $x \in HPTP^r$.  Suppose $v
\in PTP^r$.  Hence
$\O(x) \in \{PO^r(v),0\}$.
Suppose for a contradiction that $x \in PTP^r$.  Then, $PO^r(x) <
PO^r(v)$.  However, $x \overline{<} T^{PO^r(v)} x \overline{<}
T^{PO^r(x)} x, x \lessdot T^{PO^r(x)}x$ and Lemma \ref{c1} imply that
$PO^r(x) = \frac{PO^r(v)|x|}{|v|} > PO^r(v)$, a contradiction.  So,
the identity $P^r(v) + P^r(x) = 1$ easily follows (since $v \in
EPTP^r$).

The rest of the proof is almost exactly the same as the proof of the
first part of Lemma \ref{pq}, getting rid of $u$.
Again, the results follow with simple
modifications in the decreasing case.
\end{proof}

\subsection{Repeated application of $\phi$ on larger collections of $T$-pairs}
\begin{ov}  In this section, we consider the longer monomial terms that
arise in the expansion of the formula for $R_J$.  As in the previous
cases, most terms cancel; we therefore explicitly give the groupings
which cancel (Corollary \ref{top}) or almost cancel, and show that
what remains is simply $R_{\text{GGS}}$, proving
Theorem \ref{mt}.  In order to do this, we need to define the notion
of a $T$-chain, which generalizes $T$-quadruples to the objects needed
to handle the longer monomials, and prove some more combinatorial
results regarding these.
\end{ov}

\begin{defe}
A {\sl $n,m$ $TP$-chain} is a chain $x = (x_1,\ldots,x_n) \in TP_+^m
\times TP_-^{n-m}$, so that $E_{x_1} \cdots E_{x_n} \neq 0$.
Let $TPC_{n,m}$ denote the set of $n,m$ TP-chains.
Let $TPC$ denote the set of all $TP$-chains,
i.e. $TPC = \cup_{i > j} TPC_{i,j}$.
\end{defe}

\begin{defe}  A chain $(x_1, \ldots, x_n) \in TPC_{n,m}$ is said to be
{\sl outer} if\linebreak $\O(x_1) > \ldots > \O(x_m)$ and
$\O(x_{m+1}) <
\ldots < \O(x_n)$. Let $TPC^o_{n,m}$ denote the set of such
$n,m$-chains, and let $TPC^o$ denote the set of all outer chains
for any $n, m$.
\end{defe}

\begin{note}
Lemmas \ref{ail}-\ref{onp} below have obvious analogues
obtained by changing
the sign of all $T$-pairs, reversing the order of $T$-chains, and
replacing $TPC_{n,m}$ with $TPC_{n,n-m}$ in all forms.  These results
are not stated but will be referred to in the same manner as the actual
results stated.
\end{note}

\begin{lemma}   \label{ail}
$($i\,$)$ Suppose $(x_1, x_2, x_3) \in TPC_{3,3}$.  Suppose
$\psi(x_2, x_3) = (y_2, y_3)$ and $\psi(x_1,y_2) = (z_1,z_2)$
where $(x_2, x_3), (x_1,y_2) \in CTQ_s^{d}$.  Then $\O(x_1) <
\O(x_2)$ iff
$\O(z_2) <
\O(y_3)$, $\O(x_1) = \O(x_2)$ iff $\O(z_2) = \O(y_3)$, and $\O(x_1) >
\O(x_2)$ iff $\O(z_2) > \O(y_3)$.

$($ii\,$)$ Suppose $(x_1, x_2) \in TPC_{2,2}$.  Suppose
$\psi^l(x_2) = (y_2, y_3)$ and $\psi(x_1, y_2) = (z_1, z_2)$
where $x_2 \in PTP^l$ and
$(x_1, y_2) \in CTQ_s^d$.  Then $\O(x_1) ? \O(x_2)$ iff $\O(z_2) ?
\O(y_3)$ for $?$ any order relation $=, <,$ or $>$.
\end{lemma}

\begin{proof} The proof is the same for both (i) and (ii), does not mention
$x_3$, and is given in the following paragraphs.
First, we note that $\O(x_1) = \O(x_2)$ iff $x_1$ and $x_2$ reverse
orientation, which is true iff $z_2$ and $y_3$ reverse orientation,
which is true iff $\O(z_2) = \O(y_3)$.

Now, we show that $\O(x_1) < \O(x_2)$ implies $\O(z_2) < \O(y_3)$.  To
reach a contradiction, suppose that $\O(x_1) < \O(x_2)$ and $\O(z_2) >
\O(y_3)$.  By nilpotency, $x_1$ and $y_3$ preserve orientation.  Thus
$z_2$ and $x_2$ also preserve orientation.  Write $\O(y_2) = p\, \O(z_1)
+ (p-1) \O(z_2)$.

First suppose $p \geq 2$.  Then $|z_2| > \frac{1}{2} |z_1|$, so that
$T^{\O(z_1)+\O(z_2)}$ is defined on $(y_3)_h$ and therefore on
$(y_3)_l$.  However, this implies that $T^{\O(y_3)+\O(z_1)}
((z_2)_l+(y_3)_l) \overline{>} (z_2)_l + (y_3)_l \overline{>}
T^{\O(z_1)+\O(z_2)} ((z_2)_l+ (y_3)_l)$, while $T^{\O(y_3) + \O(z_1)}
((z_2)_l+(y_3)_l) \gtrdot T^{\O(z_1)+\O(z_2)} ((z_2)_l+(y_3)_l)$.
This contradicts Lemma \ref{c1}.

So $p = 1$. Then $\O(x_1) = 2 \O(z_1) + \O(z_2)$ and $\O(y_2) =
\O(z_1)$.  Since $\O(x_1) < \O(x_2)$, it must be that $\O(x_2) > 2
\O(z_1) + \O(y_3) = 2 \O(y_2) + \O(y_3)$, which implies that $|y_3| >
\frac{1}{2} |y_2|$.  So $T^{\O(y_2)+\O(y_3)} [(z_2)_l + (y_3)_l]
\overline{<} (z_2)_l + (y_3)_l$.  Now, write $\O(z_1) + \O(z_2) = q [
\O(y_2) + \O(y_3) ] + m$, for $0 \leq m < \O(y_2) + \O(y_3)$.  Since
$\O(y_2) + \O(y_3) < \O(z_1) + \O(z_2) < \O(x_2) - \O(y_2)$, it
follows that $T^{m} ((y_3)_l + (z_2)_l) \lessdot (y_3)_l + (z_2)_l$
($q = 1$) or $T^m ((y_3)_l + (z_2)_l) \overline{<} (y_3)_l + (z_2)_l$
($q > 1$).  By Lemma \ref{c1}, it follows that $c_{(y_3)_l + (z_2)_l}$
divides both $m$ and $\O(y_2) + \O(y_3)$.  But then $c_{(y_3)_l +
(z_2)_l}$ divides $\O(x_1) - \O(z_1)$ and $\O(x_2) - \O(z_1)$, which
shows that $\O(x_1) > \O(x_2)$, contrary to assumption.

Next, suppose $\O(x_1) > \O(x_2)$.  We show $\O(z_2) > \O(y_3)$.
Find $p$ and $q$ such that $\O(x_2) = (p+1) \O(y_2) + p\, \O(y_3)$ and
$\O(x_1) = (q+1) \O(z_1) + q \, \O(z_2)$.  In this case, $\O(y_2) =
q\, \O(z_1) + (q-1) \O(z_2)$.  Hence, $\O(x_2) = (p+1) q\, \O(z_1) +
(p+1)(q-1) \O(z_2) + p\, \O(y_3)$.  By assumption, $(1 - pq) \O(z_1)
+ (p - pq + 1) \O(z_2) - p\, \O(y_3) > 0$.  Since $pq \geq 1$, this
in particular implies that $\O(z_2) > \O(y_3)$, as desired.
\end{proof}

\begin{lemma}  \label{ail2} $($i\,$)$ Suppose $(x_1, x_2, x_3) \in
TPC_{3,1}$  with $|x_1| > |x_2| + |x_3|$ and \\ $\O(x_2) <
\O(x_3)$.  Let
$\psi(x_2,x_3) = (y_2,y_3)$ and suppose $\O(y_2) < \O(x_1)$.  Let
$\psi(x_1,y_2) = (z_1, z_2)$.  Then $\O(z_2) < \O(y_3)$.

$($ii\,$)$ Similarly, suppose $(x_1, x_2) \in TPC_{2,1}$ with
$|x_1| > |x_2|$ and $x_2 \in EPTP^r$.  Then set $\psi^r(x_2) =
(y_2, y_3)$ and suppose $\O(y_2) < \O(x_1)$. Setting $\psi(x_1,
y_2) = (z_1, z_2)$, it follows that $\O(z_2) < \O(y_3)$.
\end{lemma}

\begin{proof}
Again, parts (i) and (ii) have nearly the same proof, which follows.
Suppose, on the contrary, that $\O(z_2) \geq \O(y_3)$.  Clearly
$\O(z_2) \neq \O(y_3)$, else $y_3$ would have reversed orientation,
which is not possible by Lemma \ref{nc}. So $\O(z_2) > \O(y_3)$.  Now,
$T^{\O(x_1)}$ is defined on $(y_3)_l$ since, in case (i),
$T^{\O(x_1)}$ is defined on $(x_2)_h$ and $(x_3)_h$, which follows
from $|x_1| > |x_2|+|x_3|$, and in case (ii), $T^{\O(x_1)}$ is defined
on $(x_2)_h$.  In particular, $T^{\O(z_2)}$ is defined on $(y_3)_l$.
But now, $T^{-\O(z_1)}((y_3)_l+(z_2)_l) \overline{<}
T^{\O(y_3)}((y_3)_l+(z_2)_l)$ but $T^{\O(z_2)} ((y_3)_l+(z_2)_l)
\lessdot T^{\O(y_3)}((y_3)_l + (z_2)_l)$, contradicting Lemma
\ref{c1}.
\end{proof}

\begin{lemma} \label{pmol}
Suppose $x = (x_1, x_2, x_3, x_4) \in TPC^{o}_{4,2}$ with $|x_2| <
|x_3| + |x_4|$ and $|x_3| < |x_2| + |x_1|$.  Let $\psi(x_1, x_2) = (y_1,
y_2)$ and $\psi(x_3,x_4) = (y_3, y_4)$. Then $|x_2| > |x_3|$ iff $(y_1,y_2,
x_3, x_4) \in TPC^o$ and $|x_3| > |x_2|$ iff $(x_1, x_2, y_3, y_4) \in TPC^o$.
\end{lemma}

\begin{proof} First we note that, if $|x_2| > |x_3|$, then $(y_1, y_2, x_3,
x_4) \in TPC^o$, for the following reason.  Suppose $|x_2| > |x_3|$ and
write $\phi(x_2, x_3) = (z_2, z_3)$ or $\phi(x_2, x_3) = z_3$, depending
on whether $(x_2, x_3)$ is bad or good. Then $|z_2| + |z_3| < |x_4|$, so
Lemma \ref{ail2} applies and shows that $\O(y_2) < \O(x_3)$, as desired.
By symmetry, $|x_3| > |x_2|$ implies that $(x_1, x_2, y_3, y_4) \in TPC^o$.

So, it suffices to show that $|x_2| > |x_3|$ if $(x_1, x_2,
y_3, y_4) \in TPC^o$ and $|x_3| > |x_2|$ if $(y_1, y_2, x_3, x_4) \in TPC^o$.
By the symmetry of the situation, we need only prove the first.

Assume, for sake of contradiction, that $|x_3| > |x_2|$ and $(x_1,
x_2, y_3, y_4) \in TPC^o$.  First suppose that $|x_1| + |x_2| > |x_3|
+ |x_4|$.  Then we note that $T^{-\O(x_2)}$ is defined on $(y_4)_l$
because it is defined on $(x_3)_l$ and $(x_4)_l$.  Next we note that
$T^{-\O(x_2)} (y_4)_l \overline{>} T^{-\O(y_3)} (y_4)_l$ because
$|y_4| > |y_3| + |x_2|$.  Also, $T^{-\O(x_2)} (y_4)_l \not \perp (y_4)_h$,
and $T^{-\O(y_3)} (y_4)_l \overline{<} (y_4)_h$. This contradicts
Lemma \ref{c1}.

On the other hand, it is impossible that $|x_3| + |x_4| > |x_1| +
|x_2|$.  If this were true, then it would follow that $T^{\O(x_2) +
\O(x_3)}$ shifted $T^{-\O(x_2)} ((x_1)_h + (x_2)_h)$ to the right by
$|x_3| - |x_2| < |x_1|$, while $x_1 \prec T^{-\O(x_2)} ((x_1)_h)$.
This would contradict Corollary \ref{c3}.
\end{proof}

\begin{lemma} \label{bps} Suppose $(x_1, x_2) \in TPC_{2,1}$ and $|x_1| >
|x_2|$.
Suppose further that $x_1 \in PTP^l$ and $\psi^l(x_1) = (y_1, y_2)$.
Then $\O(y_2) < \O(x_2)$.
\end{lemma}

\begin{proof}  Clearly $\O(y_2) = \O(x_2)$ would imply that $y_2, x_2,$
and $x_1$ reverse orientation, which is not possible since $|x_1| >
|x_2|$.  Suppose instead that $\O(y_2) > \O(x_2)$.
We have that $T^{\O(x_1)}$, and hence $T^{\O(y_2)}$, is defined on
$(x_2)_h$, and hence $(x_2)_l$.  Since $\O(y_2) > \O(x_2)$, we find
that $T^{-\O(y_1)} [(x_2)_l + (y_2)_l] \overline{<} T^{\O(x_2)}
[(x_2)_l + (y_2)_l] \gtrdot T^{\O(y_2)} [(x_2)_l + (y_2)_l]$,
contradicting Lemma \ref{c1}.  So $\O(y_2) < \O(x_2)$. \end{proof}

\begin{lemma} \label{onp} Suppose $x = (x_1, x_2, x_3) \in TPC^o_{3,2}$ with
$|x_2| < |x_3| < |x_1| + |x_2|$.  Then
$x_2 \notin PTP^l$ and $x_3 \notin EPTP^r$.
\end{lemma}

\begin{proof}  First, note that $T^{-\O(x_2)}$ is defined on
$(x_3)_l$ since it is defined on $(x_2)_h + (x_1)_h$, and the former
is a subset of the latter on the diagram.  If $x_2 \in EPTP^l$, then
set $q = PO^l(x_2)$.  If $x_3 \in EPTP^r$ and $x_2 \notin EPTP^l$, set
$q = \O(x_2) + \O(x_3) - PO^r(x_3)$.  In either case, $T^{q}$ is
defined on $T^{-\O(x_2)} ((x_3)_l)$.  Indeed, $T^{\O(x_2) + \O(x_3)}$
is defined on $T^{-\O(x_2)} ((x_3)_l)$, and $q < \O(x_2) + \O(x_3)$.
Note that, by Lemma \ref{nc}, $x_3$ and $x_2$ must preserve
orientation.  Hence, $T^q$ and $T^{\O(x_2) + \O(x_3)}$ must preserve
orientation on $T^{-\O(x_2)}((x_3)_l)$.  Set $\omega = T^{-\O(x_2)}
((x_3)_l)$. We have $\omega \overline{>} T^{q} \omega \lessdot
T^{\O(x_2) + \O(x_3)} \omega = (x_3)_h$, contradicting Lemma
\ref{c1}. \end{proof}

Now, we are ready to define maps $\Phi, \Psi_{\pm}, \Psi^r, \Psi^l,
\Psi'$, and $\Phi'$ which parameterize terms which cancel in the
expansion of $J^{-1} R_{st} J^{21}$ (given in Corollary \ref{top}).

\begin{ntn} For any map $f$ taking two arguments, let $f^{i,j}$ be the
map $f$ applied to the $i$ and $j$-th components of some larger
$k$-tuple.  Similarly, for any map $g$ taking only one argument,
define $g^{(i)}$ to be the map $g$ applied to the $i$-th component of
a larger $k$-tuple.
\end{ntn}

\begin{defe}  Define $\Phi: TPC^o_{n,m} \rightarrow TPC$ for $1 \leq m < n$ as
follows.  For $x = (x_1, \ldots, x_n) \in TPC^o_{n,m}$, let $\Phi(x) =
\phi^{m,m+1}(x)$.
\end{defe}

\begin{defe} Define $\Psi_+: TPC^o_{n,m} \rightarrow TPC_{n,m-1}$ for
$m \geq 2$ by $\Psi_+(x) = \psi^{m-1,m}(x)$.  For $m \leq n - 2$,
define $\Psi_-: TPC^o_{n,m} \rightarrow TPC_{n,m+1}$ by $\Psi_-(x) =
\psi^{m+1,m+2}(x).$
\end{defe}

\begin{note}  Note that, unlike in the case of the map $\psi$ (i.e. the
case of $2,m$ chains), we can have chains on which both maps $\Psi_+$
and $\Psi_-$ are defined. We see, however, that in many cases (namely,
the hypotheses of Lemma \ref{pmol}), only one map will yield an outer
chain.
\end{note}

\begin{defe} When $x = (x_1,\ldots, x_n) \in TPC_{n,m}$ and $x_m \in PTP^l$,
$x$ is said to be {\sl left-passed}
and we write $x \in LPC$.
If $x = (x_1, \ldots, x_n) \in TPC_{n,m}$ and $x_{m+1} \in
PTP^r$ or $HPTP^r$, then $x$ is said to be {\sl right-passed} or
{\sl half-passed}, and $x \in
RPC$ or $HPC$, respectively.
\end{defe}

\begin{note} Note that, given $x = (x_1, \ldots, x_n) \in TPC_{n,m}$ with
$n \geq m+1 \geq 2$, $x_m \in
HPTP^l$ iff $x_{m+1} \in HPTP^r$, so in these cases the half-passed chains
not only have a half right-passed pair but also a half left-passed pair.
\end{note}
%

\begin{sloppypar}
\begin{defe} Define $\Psi^l: TPC^o_{n,m} \cap LPC \rightarrow
TPC_{n+1,m}$ by
$\Psi^l(x) =$ $(\psi^l)^{(m)} (x)$
Similarly, define $\Psi^r: TPC^o_{n,m} \cap
(RPC \cup HPC) \rightarrow TPC_{n+1,m+1}$ by $\Psi^r(x) =
(\psi^r)^{(m+1)} (x)$.
\end{defe}
\end{sloppypar}

As in the case of $\phi, \psi, \psi^*$, we have the following:

\begin{lemma} \label{PPi} Let $x, y \in TPC^o$. Then $x = \Phi(y)$ iff
$y =\Psi_{\pm}(x)$ for some sign $\pm$ or
$y = \Psi^*(x)$ for some superscript $* \in \{l,r\}$.
\end{lemma}

\begin{proof}  This follows easily from Proposition \ref{qc}. \end{proof}

\begin{defe} A chain $x = (x_1, \ldots, x_n) \in TPC_{n,m}$
is said to be {\sl negatively special} if $n \geq m+1 \geq 3$ and
$|x_{m+1}| > |x_{m}| + |x_{m-1}|$.  The set of such chains is denoted
$SPC_-$. Similarly, $x = (x_1, \ldots, x_n)$ is {\sl positively
special} (or in $SPC_+$) if $n \geq m+2 \geq 3$ and $|x_{m}| >
|x_{m+1}| + |x_{m+2}|$.  Set $SPC = SPC_+ \cup SPC_-$; the general
term will be simply {\sl special}.
\end{defe}

\begin{defe} A chain
$x = (x_1, \ldots, x_n) \in TPC_{n,m}$ is said to be {\sl negatively
reversed}
if $3 \leq m+2 \leq n$ and $x_{m}$ and $x_{m+2}$ both reverse orientation.
The set of such chains is $RC_-$.  A chain
$x = (x_1, \ldots, x_n) \in TPC_{n,m}$ is {\sl positively reversed}
(or in $RC_+$)
if $3 \leq m+1 \leq n$ and $x_{m-1}$ and $x_{m+1}$ both reverse
orientation.  As before, $RC = RC_+ \cup RC_-$ is the set of {\sl
reversed} chains.
\end{defe}

\begin{lemma} \label{rcl} If $x \in TPC^o_{n,m} \cap RC$, then $(n,m) = (3,1)$
or $(3,2)$. If $x \in TPC^o \cap RC$, then $\Phi(x), \Psi_{\pm}(x), \Psi^*(x)
\notin TPC^o$ for any $\pm, *$.  In fact, if $x \in TPC^o \cap RC_+$,
then $\Psi_-$ is not defined on $x$, and $\Psi_+(x) = (y_1, y_2, y_3)
\in TPC_{3,1}$ satisfies $\O(y_2) = \O(y_3)$.  Additionally, $\Phi(x)
= (z_1, z_2, z_3) \in TPC_{3,3}$ or $\Phi(x) = (z_1, z_2) \in
TPC_{2,2}$ where $\O(z_1) = \O(z_2)$.

Similarly, if $x \in TPC^o \cap
RC_-$, then $\Psi_+$ is not defined on $x$, and $\Psi_-(x) = (y_1,
y_2, y_3) \in TPC_{3,2}$ satisfies $\O(y_1) = \O(y_2)$.  Additionally,
$\Phi(x) =\linebreak (z_1, z_2, z_3) \in TPC_{3,0}$ or $\Phi(x) = 
(z_2, z_3) \in
TPC_{2,0}$ where $\O(z_2) = \O(z_3)$.
\end{lemma}

\begin{proof} This easily follows from Lemmas \ref{nc} and \ref{psq}. 
\end{proof}

\begin{corr} \label{rs} No outer chain is both reversed and special. \end{corr}

\begin{proof} Clear. \end{proof}

\begin{defe} We define maps $\Psi', \Phi': TPC^o \cap RC \rightarrow TPC^o$ as
follows. If $x \in RC_- \cap TPC^o$ where $\Psi_-(x) = (y_1, y_2,
y_3)$, set $\Psi'(x) = (Q(y_1, y_2), y_3) \in TPC_{2,1}$.  If $\Phi(x)
= (z_1, z_2, z_3)$, set $\Phi'(x) = (z_1, Q(z_2, z_3)) \in TPC_{2,0}$,
and if $\Phi(x) = (z_2, z_3)$, set $\Phi'(x) = (Q(z_2,z_3)) \in
TPC_{1,0}$.  Similarly, if $x \in RC_+ \cap TPC^o$ where $\Psi_+(x) =
(y_1, y_2, y_3)$, set $\Psi'(x) = (y_1, Q(y_2, y_3)) \in TPC_{2,1}$.
If $\Phi(x) = (z_1, z_2, z_3)$, set $\Phi'(x) = (Q(z_1, z_2), z_3) \in
TPC_{2,2}$, and if $\Phi(x) = (z_1, z_2)$, then set $\Phi'(x) =
(Q(z_1, z_2)) \in TPC_{1,1}$.
\end{defe}

\begin{lemma}  \label{pg1} (i) Take $x = (x_1, \ldots, x_n) \in
TPC^o_{n,m} \cap SPC_{\pm}$ for some fixed sign $\pm$.  Then
$\Psi_{\pm}(x) \in TPC^o$ or $\Phi(x) \in TPC^o$, but not both.
In the former case, $\Psi_{\pm}(x) \in SPC_\pm$ as well.

(ii) If $x = (x_1, \ldots, x_n) \in TPC^o_{n,m} \setminus SPC$,
for $(n,m) \notin \{(1,0), (1,1), (2,1)\}$, then either (a) $\Psi_+(x)
\in TPC^o$, (b) $\Psi_-(x) \in TPC^o$, or (c) neither.  In case (c),
$x \in RC$.  In case (a), $\Psi_+(x) \in SPC$, and in case (b),
$\Psi_-(x) \in SPC$.

(iii) If $x \in TPC^o_{n,m} \cap [LPC \cup HPC \cup RPC]$, then $x \in
SPC$.  Additionally, $\Psi^l(x) \in TPC^o \cap SPC_+$ if $x \in SPC_+ \cap
LPC$ and $\Psi^r(x) \in TPC^o \cap SPC_-$ if $x \in SPC_- \cap [RPC \cup HPC]$.
\end{lemma}

\begin{proof}  (i) This follows immediately from Lemma \ref{ail} in the
cases that $m \geq 2, x \in SPC_+$ and $m \leq n-2, x \in SPC_-$.
(Note that it is impossible to have $x \in RC$ by Corollary \ref{rs}.)
In the case that $m = 1$ and $x \in SPC_+$ or in the case that
$m = n-1$ and $x \in SPC_-$, it is
clear that $\Phi(x) \in TPC^o$.

(ii) It follows that only one of the $\Psi_{+}(x), \Psi_{-}(x)$ can be
in $TPC^o$ from Lemma \ref{pmol}.  If $x \in RC \cap TPC^o$, then no
$\Psi_{\pm}(x)$ is outer by Lemma \ref{rcl}.  By Lemma \ref{pmol}, if 
$n \geq 4$, there is nothing to prove. If $n = 3$ and $x \notin RC$, 
then Lemma \ref{ail2} gives the desired result.  For $n = 2$ the 
result is easy.

(iii) The two statements follow from Lemmas \ref{onp}, \ref{bps},
respectively. \end{proof}

\begin{defe} We define the following subsets of $TPC^o$: \\
(a) $C_a = \{x \in SPC \cap TPC^o \mid \Phi(x) \notin TPC^o \}$ \\
(a') $C_b = \{x \in SPC \cup TPC_{2,1} \mid x, \Phi(x) \in TPC^o \}$ \\
(b) $C_c^{\pm} = \{x \in TPC^o \setminus [SPC \cup RC] \mid
\Psi_{\pm}(x) \in TPC^o \}$\\

Furthermore, let $C_a^{\pm} = C_a \cap SPC_{\pm}$, $C_b^{\pm} = C_b
\cap SPC_{\pm}$, and $C_c = C_c^+ \cup C_c^-$.
\end{defe}

\begin{lemma} \label{ds} The subsets $C_a, C_b, C_c, RC \cap TPC^o,
TPC_{1,1},$ and
$TPC_{1,0}$ are disjoint and their union is all of $TPC^o$.
\end{lemma}

\begin{proof}
That $C_a, C_b$, and $C_c$ are disjoint from $RC$ follows from
Corollary \ref{rs}, Lemma \ref{rcl}, and the definition, respectively.
The other facts regarding disjointedness are obvious.  To check that
the union is all of $TPC^o$, we apply Lemma \ref{pg1}.ii, which shows
that $C_c = TPC^o \setminus [SPC \cup RC \cup TPC_{1,1} \cup TPC_{1,0}
\cup TPC_{2,1}]$. It is clear, though, that $TPC_{2,1} \subset C_b$.
This proves the desired result. \end{proof}

\begin{lemma} \label{ctc}
$($i\,$)$ For each choice of sign, $\Psi_{\pm}$ maps $C_a^{\pm}
\cup C_c^{\pm}$ injectively to $C_b^{\pm}$.   \\
$($ii\,$)$ $\Psi^r$ maps $[C_a^- \cup
C_b^-] \cap [RPC \cup HPC]$ injectively to $C_b^- \setminus [RPC \cup
HPC]$ and $\Psi^l$ maps $[C_a^+ \cup C_b^+] \cap LPC$ injectively to
$C_b^+ \setminus LPC$.  \\
$($iii\,$)$ $C_c \cap [LPC \cup RPC
\cup HPC] = \emptyset$.  \\
$($iv$)$ $\Psi^r[(C_a^+ \cup C_b^+) \cap (RPC
\cup HPC)]$ and $\Psi^l[(C_a^- \cup C_b^-) \cap LPC]$ are both
disjoint from $C_b$. \\
$($v$)$ $\Psi_{\pm} (C_b)$ is disjoint from
$TPC^o$ for each choice of sign.
\end{lemma}

\begin{proof}
(i)
This follows immediately from Lemmas \ref{pg1}.i and \ref{PPi}.

(ii) This follows immediately from Lemmas \ref{pg1}.ii, \ref{PPi}, and
\ref{ppp}.

(iii) This follows immediately from Lemma \ref{onp}.

(iv) This follows immediately from the definition of $SPC$.

(v) This follows from Lemma \ref{pg1}.i, since we know $C_b$ is disjoint
from $RC$. \end{proof}

\begin{corr} \label{top} $TPC^o$
is partitioned by all sets of
the form: \\
$($a$)$ For $x \in C_a^+ \cup C_c^+$,
\begin{multline} \label{ax1} A_x =\\ \{x, \Psi_+(x)\} \cup
\{ y \mid y = \Psi^l(x) \text{\ or\ } y = \Psi^l \circ \Psi_+ (x) \}
\cup \{ x_p \in RC_+ \mid \Phi'(x_p) = x \},
\end{multline}
$($a'\,$)$ Those of the form \eqref{ax1}, replacing $+$ with $-$
and
$l$ with $r$. (Call these also $A_x$, distinguished from the
above by the ``sign''  of $x$.)  \\
$($b$)$ For a fixed $x \in TPC_{1,1}$, \begin{multline}
\label{ax2} A_x =\\ \{ x \} \cup \{ y \mid y = \Psi^l(x) \text{\
or\ } y =
\Psi^l \circ \Psi^l (x) \} \cup \{x_p \in RC_+ \mid \Phi'(x_p)
= x \}
\end{multline}
$($b'\,$)$ Those of form \eqref{ax2}, replacing $+$ with $-$, $l$
with
$r$, and $TPC_{1,1}$ with $TPC_{1,0}$.
Call these $A_x$ as well.

Here it is meant that, if $\Psi^*$ is not defined on $x$, then $y \neq
\Psi^*(x)$ for any $y$.
\end{corr}

\begin{proof} For every $x \in C_a^+ \cup C_c^+$, we see from the
Lemma, part (i), that $\Psi_+(x) \in C_b \subset TPC^o$.  Furthermore,
$\Psi^l(x) \in C_b$ and $\Psi^l \circ \Psi_+(x) \in C_b$ where these
are defined (Lemma, part ii). Now take any $y \in C_b^+$.  First note
that either $\Psi_+(\Phi(y)) = y$ or $\Psi^l(\Phi(y)) = y$ because $y
\in SPC_+$.  Note that $\Phi(y) \notin RC$ because
$\Psi_{+}(RC_{+})$ and $\Psi^l(RC_+ \cap LPC)$ are disjoint
from $TPC^o$ (Lemma \ref{rcl}).
Now, if $\Phi(y) \notin C_b$, then using Lemma \ref{ds}, $y \in
A_{\Phi(y)}$ (of type a or b, depending on whether $\Psi_+(\Phi(y)) =
y$ or $\Psi^l(\Phi(y)) = y$, respectively).  If $\Phi(y) \in C_b$,
then clearly $y \in C_b^+$. In this case, $\Phi^2(y) \notin C_b$ by
the Lemma, parts (ii), (iv), and (v), using Lemma \ref{PPi}.
Furthermore, in this case $\Psi^l(\Phi^2(y)) = \Phi(y)$, again from
the Lemma, parts (iv) and (v).  Then, $y \in A_{\Phi^2(y)}$. Thus, we see that
all elements $y \in C_b^+$ fit into at least one of the sets (a) or
(b).  To see that all elements in $C_b^+$ fit into at most one, note
that the way in which $y$ appears is uniquely determined by which of
the following holds: (A) $\Phi(y) \notin C_b$, $y = \Psi_{+} \circ
\Phi(y)$; (B) $\Phi(y) \notin C_b$, $y = \Psi^l \circ \Phi(y)$; (C)
$\Phi(y) \in C_b$, $\Phi(y) = \Psi^l \circ \Phi^2(y)$. Now, it remains
to consider those elements of $RC$.  For any $z \in RC$, we know that
$\Phi'(z) \in TPC^o_{2,2} \cup TPC^o_{2,0} \cup TPC^o_{1,1} \cup
TPC^o_{1,0}$ by Lemma \ref{rcl}.  The former two are subsets of $C_c$,
so we know that $z \in A_{\Phi'(z)}$ for any $z \in RC$, and clearly
$\Phi'(z)$ is unique. The negative case is almost identical to this
one. \end{proof}

\begin{defe} For $x \in TPC^o_{n,m}$, set
$E(x) = 1$ if $x_m \in HPTP_l$ or $x_{m+1} \in HPTP_r$ and $0$
otherwise.
Now, given $x=(x_1,
\ldots, x_n) \in TPC^o_{n,m}$, set $P(x) = (-1)^m (q-q^{-1})^n
q^{E(x)} \prod_{i = 1}^n (-q)^{S_{x_i} C_{x_i} (|x_i| - 1)} q^{P^a(x_i)}$.
\end{defe}

%
\begin{lemma} $J^{-1} R_{st} J^{21} = \sum_{x \in TPC^o} P(x) E_{Q(x)}
+ \sum_{i \neq j} e_{ii} \o e_{jj} + \sum_i q e_{ii} \o e_{ii}$.
\label{jrjl}
\end{lemma}

\begin{proof}  The terms in the expansion of $J^{-1} R_{st} J^{21}$ are
\begin{equation} \label{l219e}
q^{E(x)} \prod_{i=1}^n E_{x_i} S_{x_i} (q-q^{-1}) (-q)^{-S_{x_i}
C_{x_i} (|x_i| - 1)} q^{P^a(x_i)},
\end{equation}
for each $x = (x_1, \ldots, x_n) \in TPC^o$, as well as
the terms $e_{ii} \o e_{jj}$ for $i \neq j$ and
$q e_{ii} \o e_{ii}$ for each $i$.
\end{proof}

\begin{lemma}\ \label{cop} $($i\,$)$ If $x \in C_a \cup
C_c$, then
$\sum_{y \in A_x} P(y) E_{Q(y)} = 0$. $($ii\,$)$ If $x
=\nolinebreak (x_1)\linebreak
\in TPC_{1,1} \cup TPC_{1,0}$,  then $\sum_{y \in A_x} P(y)
E_{Q(y)} = -S_{x_1} q^{S_{x_1} [  P^s(x_1) +
C_{x_1} (|x_1| - 1)} E_{x_1}$.
\end{lemma}

\begin{proof}  (i) Let $x \in C_a^+ \cup C_c^+$. The negative case is 
similar (see comments at the end of the proof of this part).
Let $x = (x_1, \ldots, x_n) \in TPC_{n,m}$ and let
\begin{multline}
F =\\ (-1)^{C_{x_{m-1}} (|x_{m-1}| - 1)}\prod_{i =
1}^{m-2} -(-q)^{C_{x_i}(|x_i|-1)} q^{P^a(x_i)} \prod_{j = m+1}^{n}
(-q)^{C_{x_j}(1 - |x_j|)} q^{P^a(x_j)}
\end{multline}
  be the part of the formula for $P(x)$ which will not change upon
applying $\Psi_+$, $\Psi^l$, and $\Psi'$ (where applicable).

Let $B = P(\Psi^l(x))$ if $x \in LPC$ and $B = 0$ otherwise.  Set
$(y_{m-1}, y_{m}) = \psi(x_{m-1}, x_{m})$.  First we show that
\begin{equation} \label{e31}
P(x) + B = q^{C_{x_{m-1}} (|x_{m-1}| - 1) + P^r(x_{m-1}) - P^l(y_{m})} F.
\end{equation}
First, by Lemma \ref{pq},
\begin{equation} \label{xmm1}
P^l(x_{m-1}) = P^l(x_m) + P^r(x_m) + P^l(y_m).
\end{equation}
Hence,
\begin{equation} \label{xma}
P^a(x_{m-1}) + P^a(x_m) = P^r(x_{m-1}) - P^l(y_m) - 2 P^l(x_m).
\end{equation}
If $x \notin LPC$, then $B = 0$ and $P^l(x_m) \in \{0, \frac{1}{2}\}$.
In this case, \eqref{e31} follows from the definition of $E(t)$.
If $x \in LPC$, then $P^l(x_{m}) = 1$, and \eqref{xmm1} implies
$P^r(x_m) = 0$.  Then Lemma \ref{ppp} shows that
$B = (q^2 - 1) P(x)$ so that the left-hand side of \eqref{e31}
is $q^{2 + P^a(x_{m-1}) + P^a(x_m) + C_{x_{m-1}} (|x_{m-1}| - 1)} F$.
By \eqref{xma}, this is the same as the right-hand side.

Next, we show that
\begin{equation} \label{e32}
\sum_{\Phi'(y) = x} P(y) = q^{P^r(x_{m-1}) - P^l(y_{m})}
[q^{C_{x_{m-1}} (1 - |x_{m-1}|)} - q^{C_{x_{m-1}} (|x_{m-1}| - 1)}] F.
\end{equation}
Naturally, we may assume that $x_{m-1}$ reverses orientation.  Since,
in this case, $y_m$ also reverses orientation, it must be that $m = n
= 2$.  Suppose now that $x_1 = (e_i - e_j, e_l - e_{i+l-j})$ and $x_2
= (e_j - e_k, e_{i+l-j} - e_{i+l-k})$.  Now, for every $p$ such that
$1 \leq p \leq j-i-1$, set $u_p = (e_i - e_{i+p}, e_l - e_{l-p})$,
$v_p = (e_{i+p} - e_{k+j-i-p}, e_{l-p} - e_{2i+p+l-j-k})$, and $w_p =
(e_{k+j-i-p} - e_k, e_{2i+p+l-j-k} - e_{i+l-k})$.  Then $t_p = (u_p, v_p,
w_p)$ are exactly those reversed chains such that $\Phi'(t_p) = x$.
Now, we consider the possible passing properties of the $T$-pairs
involved.  Note first that, since $\O(u_p) = \O(x_1)$, $P^r(x_1) =
P^r(u_p)$ for all $p$. For the same reason, $P^l(w_p) = P^l(y_2)$ for
all $p$.  Next, note that $P^l(u_p) = 0$ for all $p$---otherwise,
applying $T^{PO^l(u_p)}$ to $x_1$ would contradict nilpotency.  For the
same reason, $P^r(w_p) = 0$ for all $p$.  Next, note that $P^l(v_p) =
P^r(v_p) = 0$ for all $p$. This follows from the fact that $P^l(u_p) =
0$ for all $p$, making use of Lemma \ref{pq}.  Now, \eqref{e32}
follows readily.

Let $D = P(\Psi^l \circ \Psi_+(x))$ when $\Psi_+(x) \in LPC$ and $D = 0$
otherwise.  Finally, we show that
\begin{equation} \label{e33}
P( \Psi_+ (x) ) + D = -q^{C_{x_{m-1}} (1 - |x_{m-1}|) +
P^r(x_{m-1}) - P^l(y_{m})} F.
\end{equation}
The proof is similar to the proof two paragraphs back.  If $\Psi_+(x)
\notin LPC$, then $D = 0$.  By Lemma \ref{pq}, $P^r(x_{m-1}) =
P^r(y_m) + P^r(y_{m-1}) + P^l(y_{m-1})$.  Thus $P^a(y_{m-1}) +
P^a(y_m) + 2 P^l(y_{m-1}) = P^r(x_{m-1}) - P^l(y_m)$, proving
\eqref{e33} in this case (note that $P^l(y_{m-1}) \in
\{0,\frac{1}{2}\}$ and recall the definition of $E(t)$.)  If
$\Psi_+(x) \in LPC$, then $P^l(y_{m-1}) = 1$ and by Lemma \ref{pq} (or
the equation above), $P^r(y_{m-1}) = 0$ and $P^a(y_{m-1}) + P^a(y_m) + 2 =
P^r(x_{m-1}) - P^l(y_m)$.  If we set $\psi^l(y_{m-1}) =
(z_{m-1},z_m)$, then it follows from Lemma \ref{ppp} that neither
$z_{m-1}$ nor $z_m$ is passed, so $P^a(z_{m-1}) = P^a(z_m) = 0$.
Hence, $D = (q^2 - 1) P( \Psi_+ (x) )$, so the left-hand side
of \eqref{e33} is $q^2 P( \Psi_+ (x) ) = -q^{P^a(y_{m-1}) + P^a(y_m) +
2 + C_{y_{m}} (1 - |y_m|)} F$, which is equal to the right-hand
side by Corollary \ref{psq} and the above analysis.

Now, putting \eqref{e31}, \eqref{e32}, and \eqref{e33} together, we
get
\begin{equation}
P(x) + P(\Psi_+(x)) + B + D + \sum_{\Phi'(y) = x} P(y) = 0,
\end{equation}
as desired.

The negative case is almost the same as the above, except that the
``error term'' $E(t)$ over-corrects, but this is counteracted by the
fact that now $\Psi^r$ is defined on $HPC$ as well as $RPC$.  (The
details are omitted.)

(ii) This is almost the same as the proof above. First we note that
the result is clear when $x \notin LPC \cup RPC \cup HPC$, as
$A_x = \{x\}$ and $P(x)$ gives the desired formula (bearing in mind the
definition of $E$).  Suppose $x = (x_1) \in
TPC_{1,1} \cap LPC$. Let $B = P(\Psi^l(x))$ and
let $D = P(\Psi^l \circ \Psi^l(x))$ if $x \in LPC$ and $\Psi^l
\in LPC$, and $D = 0$ otherwise.  Let $F = -(-1)^{C_{x_1} (1 - 
|x_1|)} (q-q^{-1})$.
Set $(y_1, y_2) = \psi^l(x_1)$.  First we show that
\begin{equation} \label{e34}
B + D = [ q^{C_{x_1} (1 - |x_{1}|) + P^r(x_1) + 1} - q^{C_{x_1}
(1 - |x_1|) + P^r(x_1) - 1} ] F.
\end{equation}

Lemma \ref{ppp} shows that $P^l(y_2) = 0$ and
\begin{equation} \label{ppn}
P^r(x_1) = P^r(y_2) + P^a(y_1) + 2 P^l(y_1).
\end{equation}
First suppose $D = 0$.  In this case, \eqref{ppn} implies \eqref{e34}
easily.  If $D \neq 0$, then $P^r(x_1) = P^l(x_1) = P^l(y_1) = 1$, and
$P^r(y_1) = P^r(y_2) = 0$. Then, $D =
(q - q^{-1})^2 q^{C_{x_1} (1 - |x_1|)} F$.  Furthermore, $B =
(1 - q^{-2}) q^{C_{x_1} (1 - |x_1|)} F$.  From this \eqref{e34} follows.

Next, we claim that
\begin{equation} \label{e35}
\sum_{\Phi'(y) = x} P(y) = (q - q^{-1}) q^{P^r(x_1) - P^l(y_2)}
[q^{C_{x_{1}} (|x_1| - 1)} - q^{C_{x_1} (1 - |x_1|)}] F.
\end{equation}
This follows by exactly the same arguments as used in the second
paragraph of the proof of part (i).  But, $P^l(y_2) = 0$, so,
putting \eqref{e34} and \eqref{e35}  together, the result easily follows.
The negative/$TPC_{1,0}$ case is almost the same as this, bearing in
mind the final comments in the proof of the previous part.  \end{proof}

{\it Proof of the main theorem, {\rm \ref{mt}}.} This follows
immediately from Lemmas \ref{jrjl} and \ref{cop} and Corollary
\ref{top}, using
\eqref{frge} for $R_{\text{GGS}}$. \qed

\section{Acknowledgements}
I would like to thank Pavel Etingof for
introducing me to this problem and advising me.  I would also like to
thank the Harvard College Research Program for their support.
Finally, I am indebted to Gerstenhaber, Giaquinto, and Hodges for
valuable discussions and for sharing some unpublished results.

\appendix

\section{Proof of Proposition \ref{gp}}

In this section we explicitly compute $R_{\text{GGS}}$ for generalized
Cremmer-Gervais triples, the only triples satisfying $|\Gamma_1| + 1 =
|\Gamma|$ (omitting only one root), thereby proving Giaquinto's
formula, \eqref{gcgr}, and thus Proposition \ref{gp}.  Recall from
Example \ref{gcge} the results given in \cite{GG}: these triples are
indexed by $(n,m)$ where $n = |\Gamma| + 1$ is the dimension, $m$ is
relatively prime to $n$, $\Gamma_1 = \Gamma \setminus
\{\alpha_{n-m}\}$, $\Gamma_2 = \Gamma \setminus \{\alpha_{m}\}$, and
$T \alpha_i = \alpha_{\Res(i+m)}$, where $\Res$ give the residue mod
$n$ in $\{1,\ldots,n\}$.  Moreover when $s$ is taken to have trace
zero in the first component, it is uniquely given by $s_{ii}^{ii} = 0$
and $s^{ij}_{ij} = \frac{1}{2} - \frac{1}{n}\Res\bigl(\frac{j-i}{m})$
where $i
\neq j$. Then, the only difficulty is in computing $q^{s} \tilde a
q^{s}$, so here we use \eqref{frge} to prove \eqref{gcgr}.

Clearly we have
\begin{equation}
q^s \tilde a q^s = \sum_{x \in TP_+}
   q^{-P^s(x)+r(x)} E_{x}^{21} - q^{P^s(x)-r(x)} E_{x},
\end{equation}
where $r(e_k - e_{i+k-j}, e_i - e_j) = s_{j,i+k-j}^{j,i+k-j} +
s_{i,k}^{i,k}$, for $j < i$, since $C_x = 0$ for all $x$.  It
suffices, then, to show $P^s(x) - r(x) = \frac{2\O(x)}{n}$
for all $x$.  Take $x = (e_k - e_{i+k-j}, e_i - e_j)$.  Below we
use the notation $[\text{statement}] = 1$ if ``statement'' is true and
0 if it is false.

One sees that
\begin{multline}
s_{j,i+k-j}^{j,i+k-j} + s_{i,k}^{i,k} = 1 - \frac{1}{2}[2j =
i+k] - \frac{1}{2}[i = k] -
\\ \frac{1}{n}\text{Res}\bigl(\frac{i+k-2j}{m}\bigr) -
\frac{1}{n}\text{Res} \bigl(\frac{k-i}{m}\bigr) \\ = 1 -
\frac{1}{2}[2j = i+k] - \frac{1}{2}[i = k] - \frac{2}{n}
\text{Res}\bigl(\frac{k-j}{m}\bigr) + M_{i,j,k},
\end{multline}
where
\begin{multline}
M_{i,j,k} = [2j \neq i+k][\text{Res}\bigl(\frac{k-j}{m}\bigr) >
\text{Res}\bigl(\frac{i+k-2j}{m})] \\
- [i \neq k][\text{Res}\bigl(\frac{k-j}{m}\bigr) <
\text{Res}\bigl(\frac{k-i}{m}\bigr)].
\end{multline}
  Thus, since $\text{Res}\bigl(\frac{k-j}{m}\bigr) = \O(x)$, it
suffices to show $\frac{1}{2}[2j = i+k] + \frac{1}{2}[i = k] + P^r(x)
+ P^l(x) = 1 + M_{i,j,k}$.  Note that $i = k$ iff $x \in HPTP^r$ and
$2j = i + k$ iff $x \in HPTP^l$, and in these cases $M_{i,j,k} = 0$,
so it suffices to consider $x \notin HPTP$.  In this case we need to
show $1 + M_{i,j,k} = P^r(x) + P^l(x)$.  Now,
\begin{multline}
1+M_{i,j,k} =\\ [\text{Res}\bigl(\frac{k-j}{m}\bigr) >
\text{Res}\bigl(\frac{i+k-2j}{m})] +
[\text{Res}\bigl(\frac{k-j}{m}\bigr) >
\text{Res}\bigl(\frac{k-i}{m}\bigr)],
\end{multline}
and it is not difficult to see that $[\text{Res}\bigl(\frac{k-j}{m}\bigr) >
\text{Res}\bigl(\frac{i+k-2j}{m})] = P^l(x)$ while
$[\text{Res}\bigl(\frac{k-j}{m}\bigr) >
\text{Res}\bigl(\frac{k-i}{m}\bigr)] = P^r(x)$.  This finishes the proof.
\qed



\begin{thebibliography}{CGS}

\bibitem[BD]{BD}A. A. Belavin  and V. G. Drinfeld, {\it Triangle
equations and simple Lie algebras}, Soviet
Sci. Rev. Sect. C: { Math. Phys. Rev.} {\bf 4} (1984), 93--165.

\bibitem[ESS]{ESS}P. Etingof, T. Schedler, and O. Schiffmann,
{\it Quantization of dynamical $r$-matrices}, to appear.

\bibitem[GG]{GG} M. Gerstenhaber and A. Giaquinto, {\it Boundary
solutions of the classical Yang-Baxter equation}, { Lett. Math.
Phys.} {\bf 40} (1997), 337--353.

\bibitem[GGS]{GGS}M. Gerstenhaber and A. Giaquinto, and S. D.
Schack, {\it Construction of quantum groups from Belavin-Drinfeld
infinitesimals},  
   Israel Math. Conf. Proc. {\bf 7} (1993), 45--64.

\bibitem[GH]{GH} A. Giaquinto and T. J. Hodges, {\it Nonstandard
solutions of the Yang-Baxter equation}, { Lett.  Math. Phys.}
{\bf 44} (1998), 67--75.

\bibitem[H]{H} T. J. Hodges, {\it Nonstandard quantum groups
associated to certain Belavin-Drinfeld triples}, {\ Contemp. Math.}
{\bf 214} (1998), 63--70.

\bibitem[H2]{H2} \bysame,  {\it The Cremmer-Gervais solutions of
the Yang-Baxter equation}, Proc. Amer. Math. Soc. {\bf 127} (1999),
  1819--1826, {\tt q-alg/9712036}.

\bibitem[S]{S} T. Schedler, {\it Verification of the GGS
conjecture for
$\mathfrak{sl}(n)$, $n \leq 12$},  {\tt math.QA/9901079.}

\bibitem[S2]{S2} \bysame , {\it On the GGS conjecture},
{\tt math.QA/003079}.

\end{thebibliography}
\end{document}